\documentclass[12pt,twoside]{article}

\usepackage{amsmath}
\usepackage{amssymb}
\usepackage{amsthm}
\usepackage{amscd}
\usepackage{color}
\usepackage{accents}
\usepackage{mathrsfs}
\usepackage{url}

\usepackage{enumitem}
\usepackage{euscript}
\usepackage{multirow}
\usepackage{latexsym,amsbsy}
\usepackage[pdftex]{graphicx}
\usepackage{epsfig}
\usepackage{subdepth}

\usepackage{footnote}

\usepackage{float}
\usepackage{color}
\usepackage{relsize}
\usepackage[geometry]{ifsym}
\usepackage{pdfpages}
\usepackage[makeroom]{cancel}
\usepackage[normalem]{ulem}

\usepackage{eqnarray}

\usepackage[linktoc=all,bookmarks,hypertexnames=false,debug,pagebackref=false,colorlinks=true,
linkcolor=blue,citecolor=blue,urlcolor=blue]{hyperref}
\usepackage{bookmark}

\usepackage{tocloft}
\usepackage[multiple]{footmisc}

\newcommand{\blfootnote}[1]{%
  \begingroup
  \renewcommand{\thefootnote}{}%
  \footnote{#1}%
  \addtocounter{footnote}{-1}%
  \endgroup}

\usepackage{titlesec}

\numberwithin{equation}{section}

\theoremstyle{plain}
\newtheorem{theorem}{Theorem}[section]

\newtheorem{lemma}[theorem]{Lemma}

\newtheorem{remark}[theorem]{Remark}

\newtheorem{example}[theorem]{Example}

\makeatletter
\renewcommand*{\@fnsymbol}[1]{\@arabic{#1}}
\newcommand*\bigcdot{\mathpalette\bigcdot@{.8}}
\newcommand*\bigcdot@[2]{\mathbin{\vcenter{\hbox{\scalebox{#2}{{\hskip 1pt}$\m@th#1\bullet$}}}}}
\makeatother

\def\C {\mathbb{C}}
\newcommand{\E}{\operatorname{\mathbb{E}}}

\def\P {\operatorname{\mathbb{P}}}
\def\R {\mathbb{R}}
\def\Z {\mathbb{Z}}

\def\Cov {\mathrm{Cov}}

\newcommand{\dd}{{{\hskip2pt} \mathrm{d}}}
\def\i{\mathrm{i}\hskip1pt}

\newcommand\sgn{\operatorname{sgn}}

\newcommand{\mcalP}{\mathcal{P}}
\newcommand{\mcalQ}{\mathcal{Q}}
\newcommand{\mcalS}{\mathcal{S}}

\newcommand{\stirlingI}[2]{\genfrac{[}{]}{0pt}{}{#1}{#2}}

\newcommand{\stirlingII}[2]{\genfrac{\{}{\}}{0pt}{}{#1}{#2}}

\def\Fbar{\mathop{\bar{F}}}

\newcommand{\eqdist}{\overset{\mathcal{L}}{=}}

\begin{document}

\title{\vspace{-50pt}
\Large\textbf{Cumulants, Partitions, and Generalizations of P\'olya's Characterization of the Normal Distribution}
}

\author{{Armine Bagyan}$^{1,2}$ {and Donald Richards}$^{1,3}$}

\date{ }

\maketitle

\phantom{a}
\vspace{-35pt}

\begin{abstract}
Let $X_1$ and $X_2$ be independent, identically distributed random variables with finite variance.  P\'olya (1923) proved the celebrated result that if $\sqrt{2} X_1 \eqdist X_1 + X_2$ then $X_1$ is normally distributed.  P\'olya's theorem has been extended in many ways, e.g., to the case of random coefficients $R_1,R_2$ such that, under certain regularity conditions, if $(R_1^2 + R_2^2)^{1/2} X_1 \eqdist R_1 X_1 + R_2 X_2$, then $X_1$ is normally distributed.  

In this article, P\'olya's theorem is generalized \textit{via} the combinatorics of $\Pi_n$, the set-partitions of $\{1,\ldots,n\}$.  Let $X_j, R_j$, $j=1,\ldots,n$, be mutually independent random variables, where $X_1,\ldots,X_n$ are identically distributed.  For each $\nu \in \Pi_n$, a random variable $W_\nu$ is constructed from $\{X_j, R_j, j=1,\ldots,n\}$ and, subject to regularity conditions, the normality of $X_1$ is deduced from functional equations involving linear combinations of the characteristic functions of the $W_\nu$.  For $n = 3$, the $W_\nu$ are $(R_1^2 + R_2^2 + R_3^2)^{1/2} X_1$, $(R_1^2 + R_2^2)^{1/2} X_1 + R_3 X_3$, $(R_1^2 + R_3^2)^{1/2} X_1 + R_2 X_2$, $R_1 X_1 + (R_2^2 + R_3^2)^{1/2} X_2$, and $R_1 X_1 + R_2 X_2 + R_3 X_3$, corresponding to the partitions $\{123\}$, $\{12|3\}$, $\{13|2\}$, $\{1|23\}$, and $\{1|2|3\}$, respectively.  By modifying the regularity conditions, we also obtain characterizations of some contaminated normal distributions.  

A characterization of the exponential distribution through the minimum of a random sample leads to the consideration of mutually independent random variables $X_j \ge 0$, $R_j > 0$, $j=1,\ldots,n$, where $X_1,\ldots,X_n$ are identically distributed. For each $\nu \in \Pi_n$ a random variable $\widetilde{W}_\nu$ is constructed from $\{X_j, R_j, j=1,\ldots,n\}$ and, subject to regularity conditions, the exponentiality of $X_1$ is derived from functional equations based on linear combinations of the survival functions of the $\widetilde{W}_\nu$, $\nu \in \Pi_n$.  



\smallskip
\noindent
{\em Running head}: Cumulants, Partitions, and P\'olya's Characterization.

\blfootnote{\hspace{-24.5pt}
$^1$Department of Statistics, Pennsylvania State University, University Park, PA 16802, U.S.A. \\
$^2$E-mail address: aub171@psu.edu \\
$^3$E-mail address: dsr11@psu.edu}

\end{abstract}



\section{Introduction}
\label{sec_introduction}

Let $X_1$ and $X_2$ be independent, identically distributed random variables that have finite variance.  P\'olya \cite{Polya} proved the famous result that if $\sqrt{2} X_1 \eqdist X_1 + X_2$ then $X_1$ is normally distributed with mean zero.  P\'olya's theorem, since its appearance more than 100 years ago, is celebrated for having revealed a new property of the normal distribution and since then for motivating extensive research on characterizations of linear functions of random vectors, goodness-of-fit tests for normality, P\'olya-type theorems on abstract spaces, and connections to random-coefficient versions of the Central Limit Theorem.  We refer to the monographs by Bryc \cite{Bryc}, Kagan, \textit{et al.}~\cite{Kagan_etal}, Bobkov, \textit{et al.}~\cite{Bobkov_etal}, and Kakosyan, \textit{et al.}~\cite{Kakosyan_etal}, and to the review articles by Ghurye and Olkin \cite{Ghurye_Olkin}, Hamedani \cite{Hamedani}, and Zaitsev, \textit{et al.}~\cite{Zaitsev_etal} for further details and references to the literature.  

The initial impetus for the present article arose from a recent result \cite[Theorem 2.5]{Bagyan_Richards} in which, within the context of random modulation, we applied P\'olya's characterization to derive a uniqueness property of the multivariate normal distribution.  Motivated by a question posed by an insightful reviewer of \cite{Bagyan_Richards}, we were led to investigate generalizations of P\'olya's theorem of the form \eqref{eq_Polya_R_n}, and those inquiries led eventually to the present article. 

It is noteworthy that, subject to regularity conditions, P\'olya's theorem remains valid more generally under the hypothesis that 
\begin{equation}
\label{eq_Polya_R_n}
(R_1^2+\cdots+R_n^2)^{1/2} X_1 \eqdist R_1 X_1 + \cdots + R_n X_n,
\end{equation}
where $X_1,\ldots,X_n,R_1,\ldots,R_n$ are non-trivial mutually independent random variables, and $X_1,\ldots,X_n$ are identically distributed; see, e.g., \cite[Chapter 3]{Bryc}, \cite[p.~29]{Kakosyan_etal}, \cite{Klebanov}, and the references provided there.  We shall provide here a combinatorial perspective on P\'olya's characterization that includes \eqref{eq_Polya_R_n} as a special instance.  

To motivate such generalizations, consider the case $n=2$ of \eqref{eq_Polya_R_n}.  Denote by $\E(\cdot)$ the expectation with respect to all underlying random variables; by applying the Fourier transform and the mutual independence of $X_1$, $X_2$, $R_1$, and $R_2$, it is seen that the stochastic equation $(R_1^2 + R_2^2)^{1/2} X_1 \eqdist R_1 X_1 + R_2 X_2$ is equivalent to 
\begin{align}
\label{eq_Polya_R_2}
\E \exp\!\big(\i t(R_1^2+R_2^2)^{1/2} X_1\big) = \E \exp(\i tR_1 X_1) \cdot \E \exp(\i tR_2 X_2),
\end{align}
for all $t \in \R$.  Suppose, temporarily, that $X_1$ is symmetric; then its characteristic function, $\psi(\cdot)$, is even, hence is a function of $t^2$, so we define the function $\phi(\cdot)$ by 
\begin{equation}
\label{eq_psi_phi}
\psi(t) = \phi(t^2) = \E \exp(\i tX_1),
\end{equation}
$t \in \R$.  Then \eqref{eq_Polya_R_2} reduces to 
\begin{equation}
\label{eq_Polya_R_2_cf}
\E \phi\big(t^2 (R_1^2+R_2^2)\big) - \E \phi(t^2 R_1^2) \cdot \E \phi(t^2 R_2^2) = 0
\end{equation}
$t \in \R$.  Since $\phi(t^2 R^2) \equiv (\phi \circ \log)(e^{t^2 R^2})$ then \eqref{eq_Polya_R_2_cf} is equivalent to 
\begin{equation}
\label{eq_E_phi_log}
\E (\phi \circ \log)\big(e^{t^2 R_1^2} e^{t^2 R_2^2}\big) - \E (\phi \circ \log)(e^{t^2 R_1^2}) \cdot \E (\phi \circ \log)(e^{t^2 R_2^2}) = 0.
\end{equation}
For any positive random variable $V$, define the 
operator $\E_t$ by 
\begin{equation}
\label{eq_Et}
\E_t(V) = \E (\phi \circ \log)(V^{t^2}) = \E \phi(t^2 \log V),
\end{equation}
$t \in \R$.  Letting $U_j = \exp(R_j^2)$, $j=1,2$, and applying \eqref{eq_Et} to \eqref{eq_E_phi_log}, we obtain 
\begin{equation}
\label{eq_Polya_R_2_Et}
\E_t (U_1 U_2) - \E_t (U_1) \cdot \E_t (U_2) = 0,
\end{equation}
and an attentive reader almost surely will have noticed that the left-hand side of \eqref{eq_Polya_R_2_Et} is reminiscent of $\Cov(U_1,U_2)$, the covariance between $U_1$ and $U_2$.  Once such an observation is made, it is natural to speculate that extensions of the covariance might lead to extensions of P\'olya's theorem.  As the joint cumulant extends the covariance, we are led naturally to consider general cumulant theory for generalizations of \eqref{eq_Polya_R_2_Et}.  

Consider, for instance, the case in which $n = 3$.  Recall from McCullagh \cite{McCullagh} that the joint cumulant of a set of random variables, $U_1,U_2,U_3$, is 
\begin{multline}
\label{eq_kappa_3}
\kappa(U_1,U_2,U_3) = \E(U_1 U_2 U_3) - \E(U_1 U_2) \E(U_3) - \E(U_1 U_3) \E(U_2) \\
- \E(U_2 U_3) \E(U_1) + 2 \E(U_1) \E(U_2) \E(U_3),
\end{multline}
whenever $U_1,U_2,U_3$ have sufficiently many finite moments.  In \eqref{eq_kappa_3}, we now replace each operator $\E(\cdot)$ by $\E_t(\cdot)$, and each $U_j$ by $\exp(R_j^2)$, $j=1,2,3$; then the right-hand side of \eqref{eq_kappa_3} becomes 
\begin{multline}
\label{eq_kappa_3_E_t}
\E_t\!\big(e^{R_1^2+R_2^2+R_3^2}\big) - \E_t\!\big(e^{R_1^2+R_2^2}\big) \, \E_t\!\big(e^{R_3^2}\big) 
- \E_t\!\big(e^{R_1^2+R_3^2}\big) \, \E_t\!\big(e^{R_2^2}\big) \\
- \E_t\!\big(e^{R_1^2}\big) \, \E_t\!\big(e^{R_2^2+R_3^2}\big) + 2 \E_t\!\big(e^{R_1^2}\big) \, \E_t\!\big(e^{R_2^2}\big) \, \E_t\!\big(e^{R_3^2}\big).
\end{multline}
By applying \eqref{eq_Et} to write each $\E_t$ operator in terms of the usual expectation $\E(\cdot)$, we find that \eqref{eq_kappa_3_E_t} becomes 
\begin{multline}
\label{eq_Polya_R_3_phi_cf}
\E \Big[\phi\big(t^2 (R_1^2+R_2^2+R_3^2)\big) - \phi\big(t^2 (R_1^2+R_2^2)\big) \phi(t^2 R_3^2) 
- \phi\big(t^2 (R_1^2+R_3^2)\big) \phi(t^2 R_2^2) \\
- \phi\big(t^2 (R_2^2+R_3^2)\big) \phi(t^2 R_1^2) + 2 \phi(t^2 R_1^2) \phi(t^2 R_2^2) \phi(t^2 R_3^2)\Big].
\end{multline}
Next we use \eqref{eq_psi_phi} to write each $\phi(\cdot)$ quantity in terms of $\psi(\cdot)$, and then we express each $\psi(\cdot)$ in the form of the characteristic function of an $X_j$.  Then \eqref{eq_Polya_R_3_phi_cf} leads us to speculate that, subject to an analyticity assumption on $\psi(\cdot)$ and growth assumptions on the moments of $R_1,R_2,R_3$, if 
\begin{multline}
\label{eq_Polya_R_3_phi_cf_2}
\E \Big[\exp\!\big(\i t (R_1^2 + R_2^2 + R_3^2)^{1/2} X_1\big) - \exp\!\big(\i t ((R_1^2 + R_2^2)^{1/2} X_1 + R_3 X_3)\big) \\
- \exp\!\big(\i t ((R_1^2 + R_3^2)^{1/2} X_1 + R_2 X_2)\big) - \exp\!\big(\i t (R_1 X_1 + (R_2^2 + R_3^2)^{1/2} X_2)\big) \\
+ 2 \exp\!\big(\i t (R_1 X_1 + R_2 X_2 + R_3 X_3)\big)\Big] = 0
\end{multline}
for all $t \in \R$, then $X_1$ is normally distributed. 

Indeed, we prove that results such as \eqref{eq_Polya_R_3_phi_cf_2} hold; further, we show that the symmetry of $X_1$ is a consequence of \eqref{eq_Polya_R_3_phi_cf_2}.  Thus we characterize the normal distribution through linear combinations of the characteristic functions of 
$(R_1^2 + R_2^2 + R_3^2)^{1/2} X_1$, $(R_1^2 + R_2^2)^{1/2} X_1 + R_3 X_3$, $(R_1^2 + R_3^2)^{1/2} X_1 + R_2 X_2$, $R_1 X_1 + (R_2^2 + R_3^2)^{1/2} X_2$, and $R_1 X_1 + R_2 X_2 + R_3 X_3$, random variables that correspond respectively to the partitions $\{\{1,2,3\}\}$, $\{\{1,2\},\{3\}\}$, $\{\{1,3\},\{2\}\}$, $\{\{1\},\{2,3\}\}$, and $\{\{1\},\{2\},\{3\}\}$.  

After our work was nearly completed, we discovered the article of Avksent'ev \cite{Avksent'ev}, who also characterized the normal distribution through the stochastic equation \eqref{eq_Polya_R_n}.  In \cite{Avksent'ev}, it was assumed that all odd moments of $X_1$ are equal to $0$, whereas we deduce those moment conditions as a consequence of \eqref{eq_Polya_R_n}.  Going beyond \cite{Avksent'ev}, we also treat functional equations that are more general than those arising from \eqref{eq_Polya_R_n}.  

In the case of general $n$, we utilize the theory of joint cumulants and the combinatorial theory of set-partitions \cite{McCullagh,Stanley}.  In crucial portions of the proofs of characterization, we will require a generalization of the Hornich-Hlawka inequality \cite{Hornich} and identities from the theory of M\"obius inversion \cite{Stanley}.  

Our results are presented as follows.  In Section \ref{sec_Polya_n}, we provide preliminary details on cumulant theory and the combinatorics of $\Pi_n$, and we use the coefficients appearing in the moment expansions of the joint cumulants to write down a functional equation for $\psi(\cdot)$, the characteristic function of $X_1$.  Subject to analyticity of $\psi(\cdot)$ and the growth restrictions on the moments of $R_1,\ldots,R_n$, we prove that all odd moments of $X_1$ equal zero, thereby proving that \eqref{eq_psi_phi} holds, with $\phi(\cdot)$ having a Maclaurin expansion in a neighborhood of $0$.  We apply the functional equation for $\phi(\cdot)$ to derive recurrence relations for the coefficients of the Maclaurin series, and we specify the initial coefficients such that the recurrence relations lead to the normality of $X_1$.  

In Section \ref{sec_modified_Polya}, we develop a second class of functional equations for $\psi(\cdot)$ by utilizing a family of coefficients that have appeared in modifications of the classical cumulants \cite{Richards,Sahi}.  Here, the resulting recurrence relations for the coefficients of the Maclaurin expansion require no specification of the initial coefficients, so that the normality of $X_1$ follows entirely from the functional equations; however, the derivations in this section are more involved than those in Section \ref{sec_Polya_n}.  

In Section \ref{sec_contaminated_normal}, we obtain characterizations of some contaminated normal distributions.  These results require intricate properties of the Chebyshev polynomials of the first kind, and they underscore the fact that the resulting characterization depends crucially on the specification of the initial coefficients in the Maclaurin expansion of $\phi(\cdot)$.  

In Section \ref{sec_exponential_distn}, characterizations of the exponential distribution are considered.  Analogous to the results of Sections \ref{sec_Polya_n}--\ref{sec_contaminated_normal}, a well-known characterization based on the sample minimum is generalized using the philosophy of cumulants and set-partitions.  Finally, in Section \ref{sec_conclusions}, we provide some concluding remarks and describe some open questions raised by the results given here.

\section{Generalizations of P\'olya's characterization}
\label{sec_Polya_n}

Throughout the article, $n$ is a positive integer, $n \ge 2$.  A \textit{partition} of the set $\{1,\ldots,n\}$ is a collection $\nu = \{B_1,B_2,\ldots\}$ of non-empty, pairwise disjoint subsets of $\{1,\ldots,n\}$ such that $B_1 \cup B_2 \cup \cdots = \{1,\ldots,n\}$.  The sets $B_1,B_2,\ldots$ are called the \textit{blocks} of $\nu$.  In writing a partition, we generally use shorthand notation; e.g., we write $\{124|36|5\}$ to denote the partition of $\{1,2,3,4,5,6\}$ into the blocks $\{124\},\{36\},\{5\}$.  

The set of all partitions of $\{1,\ldots,n\}$ is denoted by $\Pi_n$.  For $\nu \in \Pi_n$ the \textit{cardinality} of $\nu$, which is denoted by $|\nu|$, is the number of blocks in $\nu$.  For each block $B \in \nu$ we denote by $|B|$ the cardinality of $B$, and we let $\min(B) := \min\{m: m \in B\}$.  

Suppose that $u = (u_1,\ldots,u_n) \in \R^n$ and $B = (m_1,m_2,\ldots)$ is a block in a partition $\nu$, where $m_1 < m_2 < \cdots$.  We will use the notations $u_B := (u_{m_1},u_{m_2},\ldots)$ and 
\begin{equation}
\label{eq_u_B}
\|u_B\| := \Big(\sum_{m \in B} u_m^2\Big)^{1/2}.
\end{equation}

Let $U_1,\ldots,U_n$ be random variables whose joint moment-generating function exists in an open neighborhood of the origin.  The \textit{joint cumulant} of $U_1,\ldots,U_n$ is 
\begin{equation}
\label{eq_cumulant}
\kappa(U_1,\ldots,U_n) = \frac{\partial^n}{\partial t_1 \cdots \partial t_n} \log \E \exp(t_1 U_1 + \cdots + t_n U_n)\Big|_{t_1=\cdots=t_n=0}.
\end{equation}
We refer to McCullagh \cite{McCullagh} for an extensive treatment of the theory of cumulants.  In particular, a classical result \cite[Eq.~(2.9)]{McCullagh} is that 
\begin{equation}
\label{eq_cumulant_sum}
\kappa(U_1,\ldots,U_n) = \sum_{\nu \in \Pi_n} (-1)^{|\nu|-1} (|\nu|-1)! \prod_{B \in \nu} \E \bigg(\prod_{m \in B} U_m\bigg).
\end{equation}

Throughout, all random variables are assumed to be \textit{non-trivial}, i.e., $P(U = 0) < 1$, equivalently, $\E \exp(\i tU) \neq 1$ for all $t \in \R$.  

In the sequel, each characteristic function $\psi(\cdot)$ is assumed to be \textit{analytic}.  That is, $\psi(\cdot)$ can be represented in a neighborhood of $0$ by a Maclaurin series; equivalently, the underlying random variable $U$ satisfies an exponential integrability condition $\E(e^{\delta |U|}) < \infty$ for some $\delta > 0$ (Chow and Teicher \cite[p.~283]{Chow_Teicher}).  Further, all moments of $U$ are finite and are determined by $\{\psi^{(k)}(0), k=0,1,2,\ldots\}$.  

A characteristic function $\psi:\R \to \C$ is said to be \textit{symmetric} if $\psi(t) = \psi(-t)$ for all $t \in \R$.  For the case in which $\psi$ is analytic, symmetry implies that $\psi(t) = \phi(t^2)$, where $\phi(\cdot)$ is analytic and has in an open neighborhood of $0$ a Maclaurin expansion, 
\begin{equation}
\label{eq_cf_X}
\phi(u) = \sum_{j=0}^\infty a_j \frac{u^j}{j!}.
\end{equation}
Note that, in \eqref{eq_cf_X}, $a_0 = \phi(0) = 1$.  Further, if $\psi(\cdot)$ is symmetric and $U$ is the underlying (non-trivial) random variable then $a_1 < 0$; indeed, 
$$
a_1 = \phi'(0) = \tfrac12 \big[\phi(t^2)\big]''\Big|_{t=0} = \tfrac12 \big[\E\exp(\i tU)\big]''\Big|_{t=0} = -\tfrac12 \E(U^2) < 0.
$$

We now state our first generalization of P\'olya's theorem.  

\begin{theorem}
\label{thm_Polya_n}
Let $X_1,\ldots,X_n$ be identically distributed random variables, each having the analytic characteristic function $\psi(\cdot)$.  Let $R = (R_1,\ldots,R_n)$ be a vector of non-trivial random variables, all of whose moments are finite, and let $X_1,\ldots,X_n,R_1,\ldots,R_n$ be mutually independent.  
\begin{enumerate}[label=(\roman*), nosep]
\item \label{gen_polya_thm_i}
If $X_1$ is normally distributed with mean zero then 
\begin{equation}
\label{eq_Polya_R_n_cf}
 \E \sum_{\nu \in \Pi_n} (-1)^{|\nu|-1} (|\nu|-1)! \exp\!\bigg(\i t \sum_{B \in \nu} \|R_B\| X_{\min(B)}\bigg) = 0, \qquad t \in \R.
\end{equation}
{\hskip-30pt}Conversely, if \eqref{eq_Polya_R_n_cf} holds then 
\item \label{gen_polya_thm_ii}
$\psi(\cdot)$ is symmetric.  
\item \label{gen_polya_thm_iii}
In \eqref{eq_cf_X}, if $a_j = a_1^j$ for all $j=2,\ldots,n-1$ then $X_1$ is normally distributed with mean $0$ and variance $-2a_1$.  
\end{enumerate}
\end{theorem}

\begin{example}
\label{ex_n_2_3}
{\rm
Suppose that $n = 2$.  Then $\Pi_2 = \{12, 1|2\}$, i.e., $\Pi_2$ contains two partitions: $\{1,2\}$, which has one block, and $\{\{1\},\{2\}\}$, which has two blocks.  Then \eqref{eq_Polya_R_n_cf} reduces to \eqref{eq_Polya_R_2}, and Theorem \ref{thm_Polya_n} reduces to the original result of P\'olya \cite{Polya}.  

For $n=3$, we have $\Pi_3 = \{123, 12|3, 13|2, 1|23, 1|2|3\}$, and \eqref{eq_Polya_R_n_cf} reduces to \eqref{eq_Polya_R_3_phi_cf_2}.  Under the stated independence assumptions on $X_j, R_j$, $j=1,2,3$, Theorem \ref{thm_Polya_n} provides that if \eqref{eq_Polya_R_3_phi_cf_2} holds then $\psi(\cdot)$ is symmetric, i.e., all odd moments of $X_1$ are equal to $0$, and if, in \eqref{eq_cf_X}, $a_2 = a_1^2$ then $X_1$ is normally distributed.  
}\end{example}

\begin{remark}
\label{rem_DCT}
{\rm
In the sequel, we differentiate expressions such as \eqref{eq_Polya_R_n_cf} repeatedly with respect to $t$.  Each such interchange of expectation and differentiation is justified by applying the Dominated Convergence Theorem.  In the case of \eqref{eq_Polya_R_n_cf}, for instance, we need to show that for any $\nu \in \Pi_n$, 
$$
\E \Bigg|\sum_{B \in \nu} \|R_B\| X_{\min(B)}\Bigg|^{2j} < \infty,
$$
for all $j=0,1,2,\ldots$.  Applying the multinomial theorem to expand the sum, and then using the mutual independence of $X_1,\ldots,X_n,R_1,\ldots,R_n$, it follows that the expectation is finite since all moments of $X_1,\ldots,X_n$ and $R_1,\ldots,R_n$ are finite.  Therefore it will be valid to interchange the expectation and any number of derivatives, and we will do so without further explanation.  
}\end{remark}

\subsection{Preliminaries}
\label{subsec_prelim_Polya_charac_proof}

The following result is a consequence of M\"obius inversion on the lattice of partitions of $\{1,\ldots,n\}$ \cite[p.~306]{Stanley}; we provide a direct proof so as to keep the article self-contained.  Throughout, we denote by $|A|$ the cardinality of any set $A \subseteq \{1,\ldots,n\}$ and by $\mcalQ_n$ the collection of all \textit{non-empty} subsets of $\{1,\ldots,n\}$.  

\begin{lemma}
\label{lem_Moebius_inversion}
Let $h: \mcalQ_n \to \R$.  Then 
\begin{equation}
\label{eq_Moebius_inversion}
\sum_{\nu \in \Pi_n} (-1)^{|\nu|-1} (|\nu|-1)! \sum_{B \in \nu} h(B) = \sum_{\substack{A \subseteq \{1,\ldots,n\} \\ A \neq \emptyset}} (-1)^{n-|A|} h(A).
\end{equation}
In particular, 
\begin{equation}
\label{eq_zero_sum_coeffs}
\sum_{\nu \in \Pi_n} (-1)^{|\nu|-1} (|\nu|-1)! = 0.
\end{equation}
\end{lemma}

\begin{proof}
For a given non-empty set $A \subseteq \{1,\ldots,n\}$, we calculate the coefficient of $h(A)$ on the left-hand side of \eqref{eq_Moebius_inversion}.  The term $h(A)$ occurs once for every partition $\nu$ that contains $A$ as a block.  Let $T = \{1,\dots,n\} \setminus A$, and set $m = |T| = n - |A|$.  Since $A$ has been designated as a block, then all remaining elements of $\{1,\dots,n\}$ belong to a block lying entirely in $T$.  Therefore, to construct a partition $\nu$ that contains $A$, we need only to choose a partition of $T$.  

Suppose that the partition of $T$ contains $r$ blocks, then $|T| = n - |A|$ and the number of partitions of $T$ into $r$ blocks is 
$$
\stirlingII{n-|A|}{r},
$$
the Stirling number of the second kind (Comtet \cite[p.~204]{Comtet}, Stanley \cite[Section 1.9]{Stanley}).  Together with the block $A$, we obtain a partition of the set $\{1,\ldots,n\}$ with $|\nu| = r+1$ blocks, so the coefficient of each such partition is 
$$
(-1)^{|\nu|-1}(|\nu|-1)! = 
(-1)^r  r!.
$$
Thus, for each $r$, the corresponding contribution to the coefficient of $h(A)$ is 
$$
(-1)^r  r! \stirlingII{n-|A|}{r}.
$$
To calculate the total coefficient of $h(A)$, we sum this quantity over all $r$.  By applying an identity for the Stirling numbers of the second kind (\cite[p.~207, Theorem B]{Comtet}, \cite[p.~82, Eq.~(1.94d)]{Stanley}), we obtain 
$$
\sum_{r=0}^{n-|A|} (-1)^r  r! \stirlingII{n-|A|}{r} = (-1)^{n-|A|},
$$
which is the coefficient of $h(A)$ on the right-hand side of \eqref{eq_Moebius_inversion}.  

We remark that \eqref{eq_zero_sum_coeffs} can be obtained by setting $U_m = 1$, $m=1,\ldots,n$, in \eqref{eq_cumulant} and applying \eqref{eq_cumulant_sum}.  Alternatively, setting $h(A) = |A|$, and noting that $\sum_{B \in \nu} |B| = n$, then the left-hand side of \eqref{eq_Moebius_inversion} equals 
\begin{equation}
\label{eq_LHS}
n \sum_{\nu \in \Pi_n} (-1)^{|\nu|-1} (|\nu|-1)!.
\end{equation}
Since there are $\binom{n}{m}$ subsets of $\{1,\ldots,n\}$ of cardinality $m$ then the right-hand side of \eqref{eq_Moebius_inversion} reduces to 
\begin{equation}
\label{eq_RHS}
\sum_{m=1}^n (-1)^{n-m} m \binom{n}{m} = \frac{\dd}{\dd t} \sum_{m=0}^n (-1)^{n-m} \binom{n}{m} t^m \Bigg|_{t=1} = (-1)^n \frac{\dd}{\dd t} (1-t)^n \Big|_{t=1} = 0,
\end{equation}
and by equating \eqref{eq_LHS} and \eqref{eq_RHS}, we obtain \eqref{eq_zero_sum_coeffs}.  
\end{proof}

Next we obtain a generalization of the Hornich-Hlawka inequality \cite{Hornich,Ressel} by following the approach in \cite{Ressel}, and we also obtain an explicit expression for the sign of the function in \eqref{eq_HH}.  

\begin{lemma}
\label{lem_Bernstein_function}
For $n \ge 2$ and a nonnegative integer $r$, define the function 
\begin{equation}
\label{eq_HH}
f_r(u_1,\ldots,u_n) = \sum_{k=1}^n (-1)^{n-k} \sum_{1 \le m_1 < \cdots < m_k \le n} (u_{m_1} + \cdots + u_{m_k})^{(2r+1)/2},
\end{equation}
$u_1,\ldots,u_n \ge 0$.  If $u_m = 0$ for any $m=1,\ldots,n$ then $f_r(u_1,\ldots,u_n) = 0$.  

Further, if $u_1,\ldots,u_n > 0$ then $f_r(u_1,\ldots,u_n) \neq 0$ and has constant sign: 
\begin{equation}
\label{eq_sgn_fr}
\sgn(f_r(u_1,\ldots,u_n)) = 
\begin{cases}
(-1)^{n-1-r}, & r \le n-2 \\
1, & r \ge n-1
\end{cases}.
\end{equation}
\end{lemma}

\begin{proof}
Let $g_r(t) = t^{(2r+1)/2}$, $t \ge 0$.  Then 
$$
f_r(u_1,\ldots,u_n) = \sum_{k=0}^n (-1)^{n-k} \sum_{1 \le m_1 < \cdots < m_k \le n} g_r(u_{m_1} + \cdots + u_{m_k}).
$$
By repeated integration, we find that 
\begin{equation}
\label{eq_fr_integral}
f_r(u_1,\ldots,u_n) = \int_0^{u_1} \cdots \int_0^{u_n} g_r^{(n)}(t_1+\cdots+t_n) \, \dd t_n \cdots dt_1.
\end{equation}
To verify that the latter integral is finite, we need only to show that $g_r^{(n)}(t_1+\cdots+t_n)$ is integrable over an open neighborhood of the origin, e.g., that for any $\epsilon > 0$, 
\begin{equation}
\label{eq_gn_integrable}
\mathop{\idotsint}_{\substack{t_1,\ldots,t_n > 0 \\ t_1+\cdots+t_n < \epsilon}} g_r^{(n)}(t_1+\cdots+t_n) \, \dd t_1 \cdots dt_n < \infty.
\end{equation}
Since $g_r^{(n)}(t) = c_{r,n} t^{r+\frac12 - n}$, where 
\begin{equation}
\label{eq_c_rn}
c_{r,n} = \prod_{m=1}^{n} (r + \tfrac32 - m).  
\end{equation}
then the integral in \eqref{eq_gn_integrable} is a constant multiple of a well-known Liouville-Dirichlet integral \cite[Section 2]{Gupta_Richards}, 
$$
\mathop{\idotsint}_{\substack{t_1,\ldots,t_n > 0 \\ t_1+\cdots+t_n < \epsilon}} (t_1+\cdots+t_n)^{r+\frac12 - n} \, \dd t_1 \cdots dt_n = \frac{1}{(n-1)!} \int_0^\epsilon t^{r - \frac12} \dd t,
$$
which, since $r \ge 0$,  is finite.  Therefore \eqref{eq_gn_integrable} is valid.  

By \eqref{eq_fr_integral}, 
\begin{equation}
\label{eq_fr_integral_repn}
f_r(u_1,\ldots,u_n) = c_{r,n} \int_0^{u_1} \cdots \int_0^{u_n} (t_1+\cdots+t_n)^{r+\frac12 - n} \, \dd t_1 \cdots dt_n,
\end{equation}
and now it follows from \eqref{eq_fr_integral_repn} that $f_r(u_1,\ldots,u_n) = 0$ if any $u_m = 0$.  

For $u_1,\ldots,u_n > 0$, it also follows from \eqref{eq_fr_integral_repn} that $f_r(u_1,\ldots,u_n) \neq 0$ and that 
$$
\sgn\big(f_r(u_1,\ldots,u_n)\big) = \sgn(c_{r,n}).
$$
By inspecting \eqref{eq_c_rn}, we obtain \eqref{eq_sgn_fr}.  
\end{proof}

For any set $A \subseteq \{1,\ldots,n\}$, define $j_A := \sum_{m \in A} j_m$.  For variables $t_1,\ldots,t_n$, define the monomial $t_A := \prod_{m \in A} t_m$.

\begin{lemma}
\label{lemma_zero_alt_sum}
Consider a sequence $a_j \in \R$, $j=0,1,2,\ldots$, where $a_0=1$.  For nonnegative integers $j_1,\ldots,j_n$, define 
\begin{equation}
\label{eq_K_sum}
K(j_1,\ldots,j_n) := \sum_{\nu \in \Pi_n} (-1)^{|\nu|-1} (|\nu|-1)! \prod_{B \in \nu} a_{j_B}.
\end{equation}
If $j_m = 0$ for any $m \in \{1,\ldots,n\}$, then $K(j_1,\ldots,j_n) = 0$.
\end{lemma}

\noindent
\textit{Proof}.  
Supposing that $j_n = 0$ define, for each $\omega \in \Pi_{n-1}$, 
$$
a_\omega = \prod_{C \in \omega} a_{j_C}.
$$
Now suppose that $|\omega| = k$, i.e., $\omega$ has $k$ blocks.  To extend $\omega$ to a partition in $\Pi_n$, we may insert $\{n\}$ into $\omega$ either as a new block or into an existing block.  

If $\{n\}$ is inserted into $\omega$ as a new block then the resulting partition, $\{\omega,\{n\}\}$ has $k+1$ blocks and the resulting summand in \eqref{eq_K_sum} is 
$$
(-1)^{|\nu|-1} (|\nu|-1)!\Big|_{\nu = \{\omega,\{n\}\}} = (-1)^{k+1-1} (k+1-1)! \, a_\omega \cdot a_{j_{\{n\}}} \equiv (-1)^k k! \, a_\omega,
$$
since $a_{j_n} = a_0 = 1$.  On the other hand, if $\{n\}$ is inserted into any of the $k$ existing blocks in $\omega$ then the resulting summand is 
$$
(-1)^{k-1} (k-1)! \, a_\omega.
$$
Therefore, the total summands from inserting $\{n\}$ into $\omega$ is 
$$
(-1)^k k! \, a_\omega + k \cdot (-1)^{k-1} (k-1)! \, a_\omega = 0.
$$
Now summing over all $\omega \in \Pi_{n-1}$, we obtain $K(j_1,\ldots,j_n) = 0$.  

Finally, if we begin with the assumption that any $j_m = 0$ then, by modifying the above argument appropriately, we obtain the same conclusion.  
$\qed$


\subsection{Proof of Theorem \ref{thm_Polya_n}}
\label{subsec_proof_Polya_n}

\noindent
\textbf{Proof of Theorem \ref{thm_Polya_n}\ref{gen_polya_thm_i}}: 
Suppose that $X_1$ is normally distributed with mean zero; then, for some $a_1 < 0$, $\E \exp(\i tX_1) = \exp(a_1 t^2)$, $t \in \R$.  Since $X_1,\ldots,X_n$ are mutually independent then, for each $\nu \in \Pi_n$, the collection $\{X_{\min(B)}: B \in \nu\}$ are mutually independent.  Therefore 
\begin{align*}
\E \exp\!\bigg(\i t\sum_{B \in \nu} \|R_B\| X_{\min(B)} \bigg) &= \prod_{B \in \nu} \E \exp(\i t \|R_B\| X_{\min(B)}) \\
&= \E \prod_{B \in \nu} \exp(a_1 t^2 \|R_B\|\} = \E \exp\!\bigg(a_1 t^2 \sum_{B \in \nu} \|R_B\|^2\bigg).
\end{align*}
Since each $\nu$ is a disjoint union of its blocks then, by \eqref{eq_u_B}, 
$$
\sum_{B \in \nu} \|R_B\|^2 = \sum_{B \in \nu} \sum_{m \in B} R_m^2 = \sum_{m=1}^n R_m^2,
$$
which does not depend on $\nu$.  Therefore for all $t \in \R$, 
\begin{multline*}
\E \sum_{\nu \in \Pi_n} (-1)^{|\nu|-1} (|\nu|-1)! \exp\!\bigg(\i t\sum_{B \in \nu} \|R_B\| X_{\min(B)} \bigg) \\
= \bigg[\sum_{\nu \in \Pi_n} (-1)^{|\nu|-1} (|\nu|-1)!\bigg] \E \exp\!\bigg(a_1 t^2 \sum_{m=1}^n R_m^2\bigg) = 0,
\end{multline*}
by \eqref{eq_zero_sum_coeffs}, so \eqref{eq_Polya_R_n_cf} holds.  
$\qed$

\bigskip


\noindent
\textbf{Proof of Theorem \ref{thm_Polya_n}\ref{gen_polya_thm_ii}}: 
On differentiating \eqref{eq_Polya_R_n_cf} $2j+1$ times with respect to $t$ and then setting $t = 0$ we obtain, for all $j=0,1,2,\ldots$, 
\begin{equation}
\label{eq_Polya_R_n_jth_deriv}
\E \sum_{\nu \in \Pi_n} (-1)^{|\nu|-1} (|\nu|-1)! \bigg(\sum_{B \in \nu} \|R_B\| X_{\min(B)}\bigg)^{2j+1} = 0,
\end{equation}

Now suppose that $j=0$.  By the mutual independence of $X_1,\ldots,X_n,R_1,\ldots,R_n$ and the identical distributions of $X_1,\ldots,X_n$, \eqref{eq_Polya_R_n_jth_deriv} reduces to 
\begin{equation}
\label{eq_odd_moment_1}
\E \Bigg[\sum_{\nu \in \Pi_n} (-1)^{|\nu|-1} (|\nu|-1)! \sum_{B \in \nu} \|R_B\|\Bigg] \cdot \E(X_1) = 0.
\end{equation}
Applying Lemma \ref{lem_Moebius_inversion} with $h(B) = \|R_B\|$, we obtain 
$$
\sum_{\nu \in \Pi_n} (-1)^{|\nu|-1} (|\nu|-1)! \sum_{B \in \nu} \|R_B\| 
= \sum_{k=1}^n (-1)^{n-k} \sum_{1 \le m_1 < \cdots < m_k \le n} (R_{m_1}^2 + \cdots + R_{m_k}^2)^{1/2}.
$$
Therefore \eqref{eq_odd_moment_1} reduces to 
\begin{equation}
\label{eq_odd_moment_1_2}
\E \Bigg[\sum_{k=1}^n (-1)^{n-k} \sum_{1 \le m_1 < \cdots < m_k \le n} (R_{m_1}^2 + \cdots + R_{m_k}^2)^{1/2}\Bigg] \cdot \E(X_1) = 0.
\end{equation}
By applying Lemma \ref{lem_Bernstein_function}, we find that the sum 
$$
\sum_{k=1}^n (-1)^{n-k} \sum_{1 \le m_1 < \cdots < m_k \le n} (r_{m_1}^2 + \cdots + r_{m_k}^2)^{1/2}
$$
is equal to $0$ if $r_m = 0$ for any $m=1,\ldots,n$.  Further, for $r_1,\ldots,r_n \neq 0$, the sum has constant, non-zero, sign given by \eqref{eq_sgn_fr}.  Since $R_1,\ldots,R_n$ are non-trivial random variables, i.e., $P(R_m \neq 0) > 0$ for all $m=1,\ldots,n$, then 
$$
\E \Bigg[\sum_{k=1}^n (-1)^{n-k} \sum_{1 \le m_1 < \cdots < m_k \le n} (R_{m_1}^2 + \cdots + R_{m_k}^2)^{1/2}\Bigg] \neq 0,
$$
so it follows from \eqref{eq_odd_moment_1_2} that $\E(X_1) = 0$.  

As inductive hypothesis, suppose that $\E(X_1^{2k+1}) = 0$ for all $k=1,\ldots,r-1$ for a fixed $r \ge 1$.  For each $\nu \in \Pi_n$, denote the corresponding blocks by $B_1,\ldots,B_{|\nu|}$.  Applying the multinomial theorem, we obtain 
\begin{align}
\label{eq_multinomial_sum}
\bigg(\sum_{B \in \nu} \|R_B\| X_{\min(B)}\bigg)^{2r+1} 
&\equiv \Bigg(\sum_{l=1}^{|\nu|} \|R_{B_l}\| X_{\min(B_l)}\Bigg)^{2r+1} \nonumber \\
&= \sum_{j_1+\cdots+j_{|\nu|} = 2r+1} \binom{2r+1}{j_1,\ldots,j_{|\nu|}} \prod_{l=1}^{|\nu|} \|R_{B_l}\|^{j_l} X_{\min(B_l)}^{j_l}.
\end{align}
By the mutual independence of $X_1,\ldots,X_n,R_1,\ldots,R_n$, the property that $X_1,\ldots,X_n$ are identically distributed, and the disjointness of the blocks, we obtain 
$$
\E \prod_{l=1}^{|\nu|} \|R_{B_l}\|^{j_l} X_{\min(B_l)}^{j_l} = \prod_{l=1}^{|\nu|} \E(\|R_{B_l}\|^{j_l}) \E(X_{\min(B_l)}^{j_l}) = \prod_{l=1}^{|\nu|} \E(\|R_{B_l}\|^{j_l}) \E(X_1^{j_l}),
$$
and therefore 
\begin{equation}
\label{eq_Xn_Rn_multinomial_sum}
\E \bigg(\sum_{B \in \nu} \|R_B\| X_{\min(B)}\bigg)^{2r+1} = \sum_{j_1+\cdots+j_{|\nu|} = 2r+1} \binom{2r+1}{j_1,\ldots,j_{|\nu|}} \prod_{l=1}^{|\nu|} \E(\|R_{B_l}\|^{j_l}) \, \E(X_1^{j_l}).
\end{equation}

For a given index $(j_1,\ldots,j_{|\nu|})$ in the sum on the right-hand side of \eqref{eq_Xn_Rn_multinomial_sum}, suppose that $\prod_{l=1}^{|\nu|} \E(X_1^{j_l}) \neq 0$.  Since $j_1+\cdots+j_{|\nu|} = 2r+1$ then at least one of $j_1,\ldots,j_{|\nu|}$ is odd.  By the inductive hypothesis, any odd $j_m \le 2r-1$ results in $\prod_{l=1}^{|\nu|} \E(X_1^{j_l}) = 0$, so it follows that $\prod_{l=1}^{|\nu|} \E(X_1^{j_l}) \neq 0$ if and only if such $j_m = 2r+1$, and then all other indices $j_i$ are equal to $0$.  Therefore $\prod_{l=1}^{|\nu|} \E(X_1^{j_l}) \neq 0$ if and only if 
$$
(j_1,\ldots,j_{|\nu|}) \in \{(2r+1,0,\ldots,0), (0,2r+1,0,\ldots,0), \ldots, (0,\ldots,0,2r+1)\}.
$$
Applying this result to \eqref{eq_Xn_Rn_multinomial_sum}, we obtain 
\begin{align*}
\E \bigg(\sum_{B \in \nu} \|R_B\| X_{\min(B)}\bigg)^{2r+1} &= \sum_{l=1}^{|\nu|} \E\big(\|R_{B_l}\|^{2r+1}\big) \E\big(X_1^{2r+1}\big) \\
&\equiv \E \bigg(\sum_{B \in \nu} \|R_B\|^{2r+1}\bigg) \E(X_1^{2r+1}),
\end{align*}
and therefore \eqref{eq_Polya_R_n_jth_deriv} reduces to 
\begin{equation}
\label{eq_zero_induction}
\E \bigg[\sum_{\nu \in \Pi_n} (-1)^{|\nu|-1} (|\nu|-1)! \sum_{B \in \nu} \|R_B\|^{2r+1}\bigg] \cdot \E(X_1^{2r+1}) = 0.
\end{equation}
By Lemma \ref{lem_Moebius_inversion} with $h(B) = \|R_B\|^{2r+1}$, 
\begin{multline}
\label{eq_Moebius_inversion_j}
\sum_{\nu \in \Pi_n} (-1)^{|\nu|-1} (|\nu|-1)! \sum_{B \in \nu} \|R_B\|^{2r+1} \\
= \sum_{k=1}^n (-1)^{n-k} \sum_{1 \le m_1 < \cdots < m_k \le n} (R_{m_1}^2 + \cdots + R_{m_k}^2)^{(2r+1)/2}.
\end{multline}
On applying the generalized Hornich-Hlawka inequality in Lemma \ref{lem_Bernstein_function}, we deduce that the sum on the right-hand side of \eqref{eq_Moebius_inversion_j} is identically zero if any $R_m = 0$.  On the other hand, if $R_1,\ldots,R_n \neq 0$ then the sum never equals zero since, by \eqref{eq_sgn_fr}, it has constant non-zero sign.  Therefore 
$$
\E \bigg[\sum_{\nu \in \Pi_n} (-1)^{|\nu|-1} (|\nu|-1)! \sum_{B \in \nu} \|R_B\|^{2r+1}\bigg] \neq 0,
$$
so it follows from \eqref{eq_zero_induction} that $\E(X_1^{2r+1}) = 0$.  Hence, by induction, $\E(X_1^{2j+1}) = 0$ for all $j=0,1,2,\ldots$.  

Since $\E(X_1^{2j+1}) = 0$ for all $j=0,1,2,\ldots$ then it follows that $\psi^{(2j+1)}(0) = 0$ for all such $j$, and therefore $\psi(t) = \phi(t^2)$, where $\phi(\cdot)$ is given in \eqref{eq_cf_X}.
$\qed$

\medskip

For the remainder of Section \ref{sec_Polya_n} and throughout Sections \ref{sec_modified_Polya}--\ref{sec_contaminated_normal}, we use the notation 
\begin{equation}
\label{eq_b_m_j}
b_{m,j} = \E(R_m^{2j})/j!,
\end{equation}
for $m=1,\ldots,n$, $j=0,1,2,\ldots$.  Note that $b_{m,j} > 0$ for all $m, j$ since $R_m$ is non-trivial.  

\bigskip


\noindent
\textbf{Proof of Theorem \ref{thm_Polya_n}\ref{gen_polya_thm_iii}}: 
Since $\psi(t) = \phi(t^2)$, then \eqref{eq_Polya_R_n_cf} reduces to 
\begin{equation}
\label{eq_Polya_R_n_cf_phi}
\E \sum_{\nu \in \Pi_n} (-1)^{|\nu|-1} (|\nu|-1)! \prod_{B \in \nu} \phi(t^2 \|R_B\|^2) = 0,
\end{equation}
$t \in \R$.  
Let $\nu \in \Pi_n$ and denote by $B_1,\ldots,B_{|\nu|}$ the blocks in $\nu$; then 
\begin{align}
\label{eq_phi_B_expansion}
\frac{\dd^{2j}}{\dd t^{2j}} \prod_{B \in \nu} \phi(t^2 \|R_B\|^2)\bigg|_{t=0} &\equiv \frac{\dd^{2j}}{\dd t^{2j}} \prod_{l=1}^{|\nu|} \phi(t^2 \|R_{B_l}\|^2)\bigg|_{t=0} \nonumber \\
&= \sum_{j_1 + \cdots + j_{|\nu|}=2j} \binom{2j}{j_1,\ldots,j_{|\nu|}} \prod_{l=1}^{|\nu|} \frac{\dd^{j_l}}{\dd t^{j_l}} \phi(t^2 \|R_{B_l}\|^2)\bigg|_{t=0}.
\end{align}
We now apply to each $\phi(\cdot)$ term in \eqref{eq_phi_B_expansion} the Maclaurin expansion \eqref{eq_cf_X}, yielding 
$$
\phi\big(t^2 \|R_B\|^2\big) = \sum_{k=0}^\infty a_k \|R_B\|^{2k} \frac{t^{2k}}{k!},
$$
for sufficiently small $t$.  Then 
$$
\frac{\dd^{j_l}}{\dd t^{j_l}} \phi(t^2 \|R_{B_l}\|^2)\bigg|_{t=0} = 
\begin{cases}
0, & j_l \textrm{ odd} \phantom{\bigg|} \\
\dfrac{j_l!}{(j_l/2)!} a_{j_l/2} \|R_{B_l}\|^{j_l}, & j_l \textrm{ even}
\end{cases}.
$$
Therefore for each $\nu \in \Pi_n$, 
\begin{align}
\label{eq_deriv_phi_B_expansion}
\frac{\dd^{2j}}{\dd t^{2j}} \prod_{B \in \nu} \phi(t^2 \|R_B\|^2)\bigg|_{t=0} &= \sum_{j_1 + \cdots + j_{|\nu|}=j} \binom{2j}{2j_1,\ldots,2j_{|\nu|}} \prod_{l=1}^{|\nu|} \frac{(2j_l)!}{j_l!} a_{j_l} \|R_{B_l}\|^{2j_l} \nonumber \\
&= (2j)! \sum_{j_1 + \cdots + j_{|\nu|}=j} \prod_{l=1}^{|\nu|} a_{j_l} \frac{\|R_{B_l}\|^{2j_l}}{j_l!}.
\end{align}
For nonnegative integers $j_1,\ldots,j_n$ and any block $B \in \nu$, define 
$$
j_B := \sum_{m \in B} j_m.
$$
By the multinomial theorem, 
\begin{equation}
\label{eq_R_B_expansion}
\frac{\|R_B\|^{2k}}{k!} = \frac{1}{k!} \bigg(\sum_{m \in B} R_m^2\bigg)^k = \sum_{j_B = k} \prod_{m \in B} \frac{R_m^{2j_m}}{m!},
\end{equation}
where the sum on the right-hand side is over all nonnegative integers $\{j_m, m \in B\}$ such that $j_B = k$.  Applying \eqref{eq_R_B_expansion} to \eqref{eq_deriv_phi_B_expansion}, multiplying over all $B \in \nu$, and collecting powers of $R_m$, we obtain 
\begin{equation}
\label{eq_deriv_R_B_expansion}
\frac{\dd^{2j}}{\dd t^{2j}} \prod_{B \in \nu} \phi(t^2 \|R_B\|^2)\bigg|_{t=0} = (2j)! \sum_{j_1+\cdots+j_n=j} \prod_{B \in \nu} a_{j_B} \cdot \prod_{m=1}^n \frac{R_m^{2j_m}}{j_m!}.
\end{equation}
Next, we substitute \eqref{eq_deriv_R_B_expansion} into \eqref{eq_Polya_R_n_cf_phi}.  Applying the mutual independence of $R_1,\ldots,R_n$ to evaluate all expectations, and using the notation \eqref{eq_b_m_j} for the moments of $R_m$, we obtain 
\begin{align*}
0 &= \E \sum_{\nu \in \Pi_n} (-1)^{|\nu|-1} (|\nu|-1)! \frac{\dd^{2j}}{\dd t^{2j}} \prod_{B \in \nu} \phi\big(t^2 \|R\|_B\big)\bigg|_{t=0} \\
&= (2j)! \sum_{j_1+\cdots+j_n=j} \bigg[\sum_{\nu \in \Pi_n} (-1)^{|\nu|-1} (|\nu|-1)! \prod_{B \in \nu} a_{j_B}\bigg] \prod_{m=1}^n b_{m,j_m},
\end{align*}
hence we obtain the recurrence relation: For $j=0,1,2,\ldots$, 
\begin{equation}
\label{eq_Polya_n_recur}
\sum_{j_1+\cdots+j_n=j} \bigg[\sum_{\nu \in \Pi_n} (-1)^{|\nu|-1} (|\nu|-1)! \prod_{B \in \nu} a_{j_B}\bigg] \prod_{m=1}^n b_{m,j_m} = 0.
\end{equation}

Suppose that $j \le n-1$; then $j_m=0$ for at least one $m$ and it follows by Lemma \ref{lemma_zero_alt_sum} that the term in square brackets in \eqref{eq_Polya_n_recur} is identically $0$.  Therefore the recurrence relation induces no restrictions on $a_1,\ldots,a_{n-1}$, and we may assign $a_1,\ldots,a_{n-1}$ freely (subject to the requirement that the resulting $\phi(\cdot)$ is a characteristic function).  Thus, we set $a_j = a_1^j$ for $j=2,\ldots,n-1$.  

For $j \ge n$, since the inner sum in \eqref{eq_Polya_n_recur} equals zero if any $j_m=0$ then the limits on the outer sum in \eqref{eq_Polya_n_recur} are restricted further to $j_1,\ldots,j_n \ge 1$.  As for the inner sum, we decompose it into the cases in which $|\nu|=1$, i.e., $\nu = \{1,\ldots,n\}$, and those for which $|\nu| \ge 2$.  If $|\nu| = 1$ then the corresponding summand is $a_{j_{\{1,\ldots,n\}}} \equiv a_{j_1+\cdots+j_n} = a_j$.  Therefore, \eqref{eq_Polya_n_recur} reduces for $j \ge n$ to 
$$
\sum_{\substack{j_1+\cdots+j_n=j \\ j_1,\ldots,j_n \ge 1}} \bigg[a_j + \sum_{\substack{\nu \in \Pi_n \\ |\nu| \ge 2}} (-1)^{|\nu|-1} (|\nu|-1)! \prod_{B \in \nu} a_{j_B}\bigg] \prod_{m=1}^n b_{m,j_m} = 0.
$$

Suppose, as inductive hypothesis, that $a_k = a_1^k$ for all $k=0,1,2,\ldots,j-1$ for a fixed $j \ge n$.  Since $j_1,\ldots,j_n \ge 1$ and $j_1+\cdots+j_n=j$ then $j_B \le j-1$ for all blocks $B \in \nu$ with $|\nu| \ge 2$.  By the inductive assumption, $a_{j_B} = a_1^{j_B}$ whenever $B \in \nu$ with $|\nu| \ge 2$, so we obtain 
$$
\prod_{B \in \nu} a_{j_B} = \prod_{B \in \nu} a_1^{j_B} = a_1^{\sum_{B \in \nu} j_B} = a_1^j.
$$
Then 
\begin{align*}
\sum_{\substack{\nu \in \Pi_n, \\ |\nu| \ge 2}} (-1)^{|\nu|-1} (|\nu|-1)! \prod_{B \in \nu} a_1^{j_B} &= \sum_{\substack{\nu \in \Pi_n \\ |\nu| \ge 2}} (-1)^{|\nu|-1} (|\nu|-1)! \: a_1^j \\
&= a_1^j \sum_{\substack{\nu \in \Pi_n \\ |\nu| \ge 2}} (-1)^{|\nu|-1} (|\nu|-1)! = - a_1^j
\end{align*}
since, by \eqref{eq_zero_sum_coeffs}, 
$$
1 + \sum_{\substack{\nu \in \Pi_n \\ |\nu| \ge 2}} (-1)^{|\nu|-1} (|\nu|-1)! \equiv \sum_{\nu \in \Pi_n} (-1)^{|\nu|-1} (|\nu|-1)! = 0.
$$
Therefore, under the inductive hypothesis, we obtain 
$$
\big(a_j - a_1^j\big) \sum_{\substack{j_1+\cdots+j_n=j \\ j_1,\ldots,j_n \ge 1}} \prod_{m=1}^n b_{m,j_m} = 0,
$$
and since $b_{m,j} > 0$ for all $m$ and $j$ then $a_j = a_1^j$.  

By induction, we obtain $a_j = a_1^j$ for all $j=0,1,2,\ldots$, and then it follows that $\phi(t^2) = \exp(a_1 t^2)$, proving that $X_1$ is normally distributed.  
$\qed$

\section{Further generalizations of P\'olya's characterization}
\label{sec_modified_Polya}

Suppose that the random variables $X_1,\ldots,X_n,R_1,\ldots,R_n$ satisfy the assumptions made in Theorem \ref{thm_Polya_n}, and consider the identity 
\begin{equation}
\label{eq_Polya_R_n_cf_mod}
\E \sum_{\nu \in \Pi_n} c_\nu \exp\!\bigg(\i t \sum_{B \in \nu} \|R_B\| X_{\min(B)}\bigg) \equiv 0,
\end{equation}
where the constants $c_\nu$ satisfy $\sum_{\nu \in \Pi_n} c_\nu = 0$.  Then it follows by a now-routine calculation that if $X_1$ is normally distributed then \eqref{eq_Polya_R_n_cf_mod} holds.  Now let us consider \eqref{eq_Polya_R_n}, the classic example of such characterizations.  

\begin{example}
\label{ex_Polya_top_bottom_parts}
{\rm 
Suppose that $c_{\{12\cdots n\}} = 1$, $c_{\{1|2|\cdots|n\}} = -1$, and $c_\nu = 0$, otherwise.  Suppose also that $X_1,\ldots,X_n,R_1,\ldots,R_n$ satisfy the same assumptions as in Theorem \ref{thm_Polya_n}.  Proceeding as in the proof of Theorem \ref{thm_Polya_n}, we deduce that all odd moments of $X_1$ are equal to $0$ and therefore $\psi(\cdot)$ is symmetric.  Further, the coefficients $a_j$ in the Maclaurin expansion of $\phi(\cdot)$ satisfy the recurrence relation 
$$
\sum_{j_1+\cdots+j_n=j} \big[a_j - a_{j_1} \cdots a_{j_n}\big] \prod_{m=1}^n b_{m,j_m} = 0,
$$
$j=0,1,2,\ldots$.  It is straightforward to prove by induction that $a_j = a_1^j$ for all $j=0,1,2,\ldots$.  Hence $\phi(t^2) = \exp(a_1 t^2)$, so $X_1$ is normally distributed.  
}\end{example}

It is noticeable that no assumptions about $a_2,\ldots,a_{n-1}$ were required in Example \ref{ex_Polya_top_bottom_parts}.  This indicates that the choice of coefficients $c_\nu$ in \eqref{eq_Polya_R_n_cf_mod} influences the number of $a_j$, $j \le n-1$, that need to be specified in order to characterize the normal distribution.  This raises the issue of whether it is possible to obtain a characterization that starts from \eqref{eq_Polya_R_n_cf_mod}, is such that $c_\nu \neq 0$ for all $\nu \in \Pi_n$, but for which no assumptions need to be made about the coefficients $a_j$, $j \le n-1$.  

In the following result, we derive such a characterization of normality.  Although the result requires no assumptions on the coefficients $a_j$, $j \le n-1$, the combinatorial details underpinning the proof are more intricate than those encountered in Section \ref{sec_Polya_n}.  

In stating the next result, we recall the notation $|B|$ for the cardinality of a block $B$ in a partition $\nu \in \Pi_n$.

\begin{theorem}
\label{thm_Polya_n_conjugate}
Let $X_1,\ldots,X_n$ be mutually independent, identically distributed random variables whose characteristic function, $\psi(\cdot)$, is analytic.  Let $R_1,\ldots,R_n$ be mutually independent random variables, all of whose moments are finite and such that $R = (R_1,\ldots,R_n)$ is independent of $X_1,\ldots,X_n$.  Then $X_1$ is normally distributed with mean zero if and only if 
\begin{equation}
\label{eq_Polya_R_n_cf_conj}
 \E \sum_{\nu \in \Pi_n} (-1)^{|\nu|-1} \bigg[\prod_{B \in \nu} (|B|-1)!\bigg] \: \exp\!\bigg(\i t \sum_{B \in \nu} \|R_B\| X_{\min(B)}\bigg) = 0, \qquad t \in \R.
\end{equation}
\end{theorem}

\medskip

\begin{example}
\label{ex_Polya_R_3_cf_conj}
{\rm
Suppose that $n = 3$.  Then Theorem \ref{thm_Polya_n_conjugate} states that, under the stated hypotheses on $X_j,R_j$, $j=1,2,3$, 
$X_1$ is normally distributed if and only if, for all $t \in \R$, 
\begin{multline}
\label{eq_Polya_R_3_cf_conj}
\E \Big[2 \exp\!\big(\i t (R_1^2 + R_2^2 + R_3^2)^{1/2} X_1\big) - \exp\!\big(\i t ((R_1^2 + R_2^2)^{1/2} X_1 + R_3 X_3)\big) \\
- \exp\!\big(\i t ((R_1^2 + R_3^2)^{1/2} X_1 + R_2 X_2)\big) - \exp\!\big(\i t (R_1 X_1 + (R_2^2 + R_3^2)^{1/2} X_2)\big) \\
+ \exp\!\big(\i t (R_1 X_1 + R_2 X_2 + R_3 X_3)\big)\Big] = 0.
\end{multline}
The functional equation \eqref{eq_Polya_R_3_cf_conj}, although similar to \eqref{eq_Polya_R_3_phi_cf_2}, results in a characterization that requires no assumption about $a_2$.  
}\end{example}

\subsection{Preliminaries}  
\label{subsec_prelim_conjugate_proof}

\begin{lemma}
\label{lem_L_n_reduction}
For $n \ge 2$ and nonnegative integers $j_1,\ldots,j_n$, define 
\begin{equation}
\label{eq_L_n_sum}
L_n(j_1,\ldots,j_n) = \sum_{\nu \in \Pi_n} (-1)^{|\nu|-1} \prod_{B \in \nu} (|B|-1)! \: a_{j_B}.
\end{equation}
Then $L_n(j_1,\ldots,j_n)$ is symmetric in $j_1,\ldots,j_n$.  

Further, if exactly $q$ of $j_1,\ldots,j_n$ are non-zero, say $j_{m_1},\ldots,j_{m_q} > 0$, then 
\begin{equation}
\label{eq_L_n_reduction_n_q}
L_n(j_1,\ldots,j_n) = 
\begin{cases}
0, & q = 0, 1 \\
\displaystyle{\frac{(n-2)!}{(q-2)!}} L_q(j_{m_1},\ldots,j_{m_q}), & q \ge 2
\end{cases}.
\end{equation}
\end{lemma}

\noindent
\textit{Proof}.  Denote by $\mcalS_n$ the group of permutations on $n$ symbols.  Also, for $\sigma \in \mcalS_n$ and $\nu \in \Pi_n$, define $\sigma(\nu) = \{\sigma(B): B \in \nu\}$, where $\sigma(B) = \{\sigma(m): m \in B\}$ for any block $B$.  The mapping $\nu \mapsto \sigma(\nu)$ is a bijection from $\Pi_n$ onto itself.  Moreover $|\sigma(\nu)| = |\nu|$, and $|\sigma(B)| = |B|$ for all $B \in \nu$.  Further, 
$$
\sum_{m \in B} j_{\sigma(m)} = \sum_{\sigma^{-1}(m) \in B} j_m = \sum_{m \in \sigma(B)} j_m = j_{\sigma(B)}. 
$$
Therefore, in the sum \eqref{eq_L_n_sum}, the summand for $\nu$ evaluated at $(j_{\sigma(1)},\ldots,j_{\sigma(n)})$ equals the summand for $\sigma(\nu)$ evaluated at $(j_1,\ldots,j_n)$.  Summing over all partitions $\nu$, we obtain 
$$
L_n(j_{\sigma(1)},\ldots,j_{\sigma(n)}) = L_n(j_1,\ldots,j_n),
$$
and since $\sigma$ was chosen arbitrarily then it follows that $L_n(j_1,\ldots,j_n)$ is symmetric in $j_1,\ldots,j_n$.  

Next, it is elementary that $L_n(j_1,\ldots,j_n) = 0$ if $q = 0,1$, so we consider the case $q \ge 2$.  Since $L_n(j_1,\ldots,j_n)$ is symmetric in $j_1,\ldots,j_n$ then we assume, without loss of generality, that $j_1,\ldots,j_q > 0$.  If $q = n$ then there is nothing to prove, so we consider the case $q = n-1$ and we will show that 
\begin{equation}
\label{eq_L_n_reduction}
L_n(j_1,\ldots,j_{n-1},0) = (n-2) L_{n-1}(j_1,\ldots,j_{n-1}).
\end{equation}
Thus suppose we choose a particular partition $\omega \in \Pi_{n-1}$, where $\omega$ has blocks $C_1,C_2,\ldots$.  Starting with $\omega$, we can construct a partition in $\Pi_n$ either by inserting $\{n\}$ as a separate block in $\omega$, or by inserting $\{n\}$ into an existing block.  If $\{n\}$ is added as a separate block in $\omega$ then we obtain a new partition, $\nu = \{\omega,\{n\}\}$, such that 
\begin{equation}
\label{eq_n_separate_block}
(-1)^{|\nu|-1} \prod_{B \in \nu} (|B|-1)! \: a_{j_B} = (-1)^{|\omega|} \prod_{C \in \omega} (|C|-1)! \: a_{j_C}.
\end{equation}
If $\{n\}$ is inserted into an existing block then we obtain $|\omega|$ possible new partitions $\nu$.  Moreover, if $\{n\}$ is inserted into the $l$th block, $C_l$, then the resulting partition $\nu$ satisfies 
$$
(-1)^{|\nu|-1} \prod_{B \in \nu} (|B|-1)! \: a_{j_B} = (-1)^{|\omega|-1} |C_l| \prod_{C \in \omega} (|C|-1)! \: a_{j_C},
$$
and by summing this over all $l$ blocks, we find that the effect of inserting $\{n\}$ into an existing block is 
\begin{equation}
\label{eq_n_same_block}
(-1)^{|\omega|-1} \Big(\sum_{C \in \omega} |C_l|\Big) \prod_{C \in \omega} (|C|-1)! \: a_{j_C} = (-1)^{|\omega|-1} (n-1) \prod_{C \in \omega} (|C|-1)! \: a_{j_C}.
\end{equation}
Therefore the total contribution (summand) associated with $\omega$ is the sum of \eqref{eq_n_separate_block} and \eqref{eq_n_same_block}, viz., 
\begin{multline}
\label{eq_total_contrib_omega}
(-1)^{|\omega|} \prod_{C \in \omega} (|C|-1)! \: a_{j_C} + (n-1) (-1)^{|\omega|-1} \prod_{C \in \omega} (|C|-1)! \: a_{j_C} \\
= (n-2) (-1)^{|\omega|-1} \prod_{C \in \omega} (|C|-1)! \: a_{j_C}.
\end{multline}
Since $\omega$ was chosen arbitrarily then, by summing \eqref{eq_total_contrib_omega} over all $\omega \in \Pi_{n-1}$, we obtain \eqref{eq_L_n_reduction}.  

Further, by iterating \eqref{eq_L_n_reduction} we obtain, for $j_1,\ldots,j_q > 0$, where $q \ge 2$, 
\begin{align*}
L_n(j_1,\ldots,j_q,0,\ldots,0) 
&= (n-2)(n-3) \cdots (q-1) L_q(j_1,\ldots,j_q) \\
&= \frac{(n-2)!}{(q-2)!} L_q(j_1,\ldots,j_q),
\end{align*}
as stated in \eqref{eq_L_n_reduction_n_q}.  
$\qed$

\smallskip

For $W = (w_1,\ldots,w_n) \in \R^n$ and $A \in \mcalQ_n$, the set of all non-empty subsets of $\{1,\ldots,n\}$, define $W_A = (w_m, m \in A)$.  Then we have the following result.

\begin{lemma}
\label{lem_R_B_conj_sign}
Let $n \ge 2$, and let $h: \mcalQ_n \to \R$.  Then 
\begin{multline}
\label{eq_h_B_sum_conj}
\sum_{\nu \in \Pi_n} (-1)^{|\nu|-1} \bigg[\prod_{B \in \nu} (|B|-1)!\bigg] \bigg(\sum_{B \in \nu} h(B)\bigg) \\
= (n-2)! \bigg[(n-1) h(\{1,\ldots,n\}) - \sum_{i=1}^n h(\{1,\ldots,n\}\setminus\{i\})\bigg].
\end{multline}
For $a > 0$, 
\begin{equation}
\label{eq_W_B_conj_sum}
\sum_{\nu \in \Pi_n} (-1)^{|\nu|-1} \bigg[\prod_{B \in \nu} (|B|-1)!\bigg] \sum_{B \in \nu} \|W_B\|^{2a} \ 
\begin{cases}
\ge 0, & a > 1 \\
  = 0, & a = 1 \\
\le 0, & 0 < a < 1
\end{cases},
\end{equation}
and for $a \neq 1$, the inequalities in \eqref{eq_W_B_conj_sum} are strict whenever at least two of $w_1,\ldots,w_n$ are non-zero.  Further, 
\begin{equation}
\label{eq_zero_sum_coeffs_conj}
\sum_{\nu \in \Pi_n} (-1)^{|\nu|-1} \prod_{B \in \nu} (|B|-1)! = 0.
\end{equation}
\end{lemma}

\noindent
\textit{Proof}.  
Fix a set $A \in \mcalQ_n$, and introduce the notation 
$$
s = |A|, \quad A^c = \{1,\ldots,n\} \setminus A, \quad m = |A^c| = n - s.
$$
To find the coefficient of $h(A)$ in the sum on the left-hand side of \eqref{eq_h_B_sum_conj}, we will total those terms in the sum that correspond to partitions $\nu$ for which $A$ is a block in $\nu$.  

Note that every such partition is of the form $\nu = \{A\} \cup \omega$, where $\omega$ is a partition of $A^c$.  Thus suppose that $|\omega| = k$; then $|\nu| = k+1$, so $(-1)^{|\nu|-1} = (-1)^k = (-1)^{|\omega|}$, and 
$$
\prod_{B \in \nu} (|B|-1)! = (s-1)! \prod_{C \in \omega} (|C|-1)!.
$$
Hence the coefficient of $h(A)$ equals 
$$
(s-1)! \sum_{\omega \in \Pi_m} (-1)^{|\omega|} \prod_{C \in \omega} (|C|-1)!.
$$
By grouping the partitions $\omega$ according to their cardinalities $k$ (i.e., according to the number of blocks in $\omega$), we find that the coefficient of $h(A)$ equals 
$$
(s-1)! \sum_{k=0}^m (-1)^k \sum_{\substack{\omega \in \Pi_m \\ |\omega| = k}} \prod_{C \in \omega} (|C|-1)!.
$$
It is well-known (\cite[p.~234, Theorem D]{Comtet}, \cite[Section 1.3]{Stanley}) that 
\begin{equation}
\label{eq_signless_Stirling}
\sum_{\substack{\omega \in \Pi_m \\ |\omega| = k}} \prod_{C \in \omega} (|C|-1)! = \stirlingI{m}{k},
\end{equation}
the signless Stirling number of the first kind.  
Indeed, for each $\omega \in \Pi_m$, the number of permutations whose cycles form the blocks of $\omega$ is $\prod_{C \in \omega} (|C|-1)!$, so the total number of permutations that have exactly $k$ cycles in their cyclic decomposition is the sum on the left-hand side of \eqref{eq_signless_Stirling}.  On the other hand $\stirlingI{m}{k}$, by definition, is the number of permutations, of the set $\{1,\ldots,m\}$, that have exactly $k$ cycles in their cyclic decomposition; therefore \eqref{eq_signless_Stirling} holds.  

Therefore the coefficient of $h(A)$ equals 
$$
(s-1)! \sum_{k=0}^m (-1)^k \stirlingI{m}{k}.
$$
A well-known identity for the signless Stirling numbers is (\cite[p.~213, Theorem A]{Comtet}, \cite[Proposition 1.3.7]{Stanley}) 
\begin{equation}
\label{eq_identity_signless_Stirling}
\sum_{k=0}^m (-1)^k \stirlingI{m}{k} \ = 
\begin{cases}
\phantom{-}1, & m = 0 \\
-1, & m = 1 \\
\phantom{-}0, & m \ge 2 
\end{cases}.
\end{equation}
Therefore only the terms corresponding to $m=0, 1$, equivalently to $s=n, n-1$, appear in the sum on the left-hand side of \eqref{eq_h_B_sum_conj}; these correspond to the full set $\{1,\ldots,n\}$ and to $\{1,\ldots,n\}\setminus\{i\}$, $i=1,\ldots,n$, respectively.  Hence the coefficient of $h(A)$ equals 
$$
(n-1)! \: h(\{1,\ldots,n\}) - (n-2)! \sum_{i=1}^n h(\{1,\ldots,n\}\setminus\{i\}),
$$
which equals the right-hand side of \eqref{eq_h_B_sum_conj}.  

In proving \eqref{eq_W_B_conj_sum}, we assume that $W \neq 0$ since the result is trivial if $W = 0$.  Setting $h(B) = \|W_B\|^{2a}$ on the left-hand side of \eqref{eq_h_B_sum_conj}, we find that the right-hand side equals 
\begin{multline}
\label{eq_W_B_conj_sum_identity}
(n-2)! \bigg[(n-1) \|W\|^{2a} - \sum_{i=1}^n \|W_{\{1,\ldots,n\} \setminus \{i\}}\|^{2a}\bigg] \\
= (n-2)! \: \|W\|^{2a} \bigg[n-1 - \sum_{i=1}^n \frac{\|W_{\{1,\ldots,n\} \setminus \{i\}}\|^{2a}}{\|W\|^a}\bigg].
\end{multline}
Let $\Sigma_n = \{(\rho_1,\ldots,\rho_n): \rho_1,\ldots,\rho_n \ge 0, \ \rho_1 + \cdots + \rho_n = 1\}$ denote the unit simplex in $\R^n$.  Setting $p_i = w_i^2 / \|W\|^2$, $i=1,\ldots,n$, then $(p_1,\ldots,p_n) \in \Sigma_n$ and, by \eqref{eq_W_B_conj_sum_identity}, it is necessary to determine the sign of 
$$
n-1 - \sum_{i=1}^n \frac{\|W_{\{1,\ldots,n\} \setminus \{i\}}\|^a}{\|W\|^{2a}} \equiv n-1 - \sum_{i=1}^n \frac{(\|W\|^2 - w_i^2)^a}{\|W\|^{2a}} = n-1 - \sum_{i=1}^n (1-p_i)^{a}.
$$

Suppose that $a > 1$.  Then $(1-p_i)^a \le 1-p_i$ for all $p_i \in [0,1]$, and the inequality is strict if $p_i \in (0,1)$.  Therefore 
$$
n - 1 - \sum_{i=1}^n (1-p_i)^a \ge n - 1 - \sum_{i=1}^n (1-p_i) = 0,
$$
with strict inequality if any $p_i \in (0,1)$.  Since at least one $p_i$ lies strictly between $0$ and $1$ whenever at least two components of $W$ are non-zero, then we have proved the first part of \eqref{eq_W_B_conj_sum}.  

For $0 < a < 1$ an argument similar to the preceding applies, with the inequalities reversed, 
and this establishes the third part of \eqref{eq_W_B_conj_sum}.  

As for the case in which $a = 1$, observe that the sum in \eqref{eq_W_B_conj_sum} is a continuous function of $a$.  Since the sum is non-positive for $a \in (0,1)$ and nonnegative for $a > 1$ then, by continuity, it equals $0$ at $a = 1$.  This proves the remaining statement in \eqref{eq_W_B_conj_sum}.  

Finally, let $W$ be a vector with at least two non-zero components.  Since $\sum_{B \in \nu} \|W_B\|^2 = \|W\|^2 \neq 0$ then, for $a = 1$, \eqref{eq_W_B_conj_sum} reduces to 
$$
\|W\|^2 \sum_{\nu \in \Pi_n} (-1)^{|\nu|-1} \prod_{B \in \nu} (|B|-1)! = 0,
$$
which establishes \eqref{eq_zero_sum_coeffs_conj}.  
$\qed$

\subsection{Proof of Theorem \ref{thm_Polya_n_conjugate}}
\label{subsec_main_thm_conj_i}

Suppose that $X_1$ is normally distributed with mean $0$.  Since \eqref{eq_zero_sum_coeffs_conj} has been proved then, by a remark made at \eqref{eq_Polya_R_n_cf_mod} \textit{infra}, we obtain \eqref{eq_Polya_R_n_cf_conj}.  

Conversely, suppose that \eqref{eq_Polya_R_n_cf_conj} holds.  Then we differentiate that identity $2j+1$ times with respect to $t$ and evaluate the outcome at $t = 0$, obtaining for all $j = 0,1,2,\ldots$, 
\begin{equation}
\label{eq_Polya_R_n_jth_deriv_conj}
\E \sum_{\nu \in \Pi_n} (-1)^{|\nu|-1} \bigg[\prod_{B \in \nu} (|B|-1)!\bigg] \bigg(\sum_{B \in \nu} \|R_B\| X_{\min(B)}\bigg)^{2j+1} = 0,
\end{equation}
Consider the case in which $j=0$.  Since $X_1,\ldots,X_n$ are mutually independent, identically distributed, and independent of $R_1,\ldots,R_n$, then \eqref{eq_Polya_R_n_jth_deriv} reduces to 
\begin{equation}
\label{eq_odd_moment_1_conj}
\E \Bigg(\sum_{\nu \in \Pi_n} (-1)^{|\nu|-1} \bigg[\prod_{B \in \nu} (|B|-1)!\bigg] \sum_{B \in \nu} \|R_B\|\Bigg) \cdot \E(X_1) = 0.
\end{equation}
By \eqref{eq_W_B_conj_sum}, 
$$
\sum_{\nu \in \Pi_n} (-1)^{|\nu|-1} \bigg[\prod_{B \in \nu} (|B|-1)!\bigg] \sum_{B \in \nu} \|R_B\| \ge 0
$$
for all $R_m \ge 0$, and the inequality is strict if at least two $R_m$ are non-zero.  Therefore 
$$
\E \sum_{\nu \in \Pi_n} (-1)^{|\nu|-1} \bigg[\prod_{B \in \nu} (|B|-1)!\bigg] \sum_{B \in \nu} \|R_B\| > 0,
$$
so it follows by \eqref{eq_odd_moment_1_conj} that $\E(X_1) = 0$.  

As inductive hypothesis, suppose that $\E(X_1^{2k+1}) = 0$ for all $k=1,\ldots,r-1$, where $r$ is a fixed integer, $r \ge 2$.  For each $\nu \in \Pi_n$, suppose that the corresponding blocks are $B_1,\ldots,B_{|\nu|}$.  Applying the multinomial theorem, we obtain 
\begin{align*}
\bigg(\sum_{B \in \nu} \|R_B\| X_{\min(B)}\bigg)^{2r+1} 
&\equiv \Bigg(\sum_{l=1}^{|\nu|} \|R_{B_l}\| X_{\min(B_l)}\Bigg)^{2r+1} \\
&= \sum_{j_1+\cdots+j_{|\nu|}=2r+1} \binom{2r+1}{j_1,\ldots,j_{|\nu|}} \prod_{l=1}^{|\nu|} \|R_{B_l}\|^{j_l} X_{\min(B_l)}^{j_l}.
\end{align*} 
By the same inductive argument that led from \eqref{eq_multinomial_sum} to \eqref{eq_Xn_Rn_multinomial_sum}, we find that \eqref{eq_Polya_R_n_jth_deriv_conj} reduces to 
\begin{equation}
\label{eq_R_B_X_1_odd_moment}
\E \bigg[\sum_{\nu \in \Pi_n} (-1)^{|\nu|-1} \bigg[\prod_{B \in \nu} (|B|-1)!\bigg] \sum_{B \in \nu} \|R_B\|^{2r+1}\bigg] \E(X_1^{2r+1}) = 0.
\end{equation}
Applying \eqref{eq_W_B_conj_sum} with $a = r + \frac12$, we obtain 
$$
\sum_{\nu \in \Pi_n} (-1)^{|\nu|-1} \bigg[\prod_{B \in \nu} (|B|-1)!\bigg] \sum_{B \in \nu} \|R_B\|^{2r+1} \le 0,
$$
and the inequality is strict whenever at least two $R_m$ are non-zero.  Since $R_1,\ldots,R_m$ are non-trivial random variables then it follows that 
$$
\E \bigg[\sum_{\nu \in \Pi_n} (-1)^{|\nu|-1} \bigg[\prod_{B \in \nu} (|B|-1)!\bigg] \sum_{B \in \nu} \|R_B\|^{2r+1}\bigg] < 0,
$$
and now applying \eqref{eq_R_B_X_1_odd_moment} we obtain $\E(X_1^{2r+1}) = 0$.  By induction, we obtain $\E(X_1^{2j+1}) = 0$ for all $j=0,1,2,\ldots$, hence the characteristic function $\psi(\cdot)$ is symmetric.  

Finally, it remains to formulate and solve the recurrence relations for the coefficients $a_j$.  Since $\psi(t) = \phi(t^2)$, then \eqref{eq_Polya_R_n_cf_conj} reduces to 
\begin{equation}
\label{eq_Polya_R_n_cf_phi_conj}
\E \sum_{\nu \in \Pi_n} (-1)^{|\nu|-1} \prod_{B \in \nu} (|B|-1)! \phi(t^2 \|R_B\|^2) = 0,
\end{equation}
$t \in \R$.  Applying the Maclaurin expansion \eqref{eq_cf_X} to each $\phi(\cdot)$ term in \eqref{eq_Polya_R_n_cf_phi_conj} and proceeding similarly to \eqref{eq_Polya_R_n_cf_phi}-\eqref{eq_Polya_n_recur}, we obtain the recurrence relation 
\begin{equation}
\label{eq_Polya_n_recur_conj}
\sum_{j_1+\cdots+j_n=j} \bigg[\sum_{\nu \in \Pi_n} (-1)^{|\nu|-1} \prod_{B \in \nu} (|B|-1)! \: a_{j_B}\bigg] \prod_{m=1}^n b_{m,j_m} = 0,
\end{equation}
$j=0,1,2,\ldots$.  

For $j=0$, \eqref{eq_Polya_n_recur_conj} reduces to 
$$
a_0 \sum_{\nu \in \Pi_n} (-1)^{|\nu|-1} \prod_{B \in \nu} (|B|-1)! = 0,
$$
which, by \eqref{eq_zero_sum_coeffs_conj}, is a tautology.  

For $j=1$, exactly one $j_r$ equals $1$ and all others are equal to $0$, so the recurrence relation \eqref{eq_Polya_n_recur_conj} reduces to 
$$
a_1 \bigg[\sum_{\nu \in \Pi_n} (-1)^{|\nu|-1} \prod_{B \in \nu} (|B|-1)!\bigg] \bigg(\sum_{m=1}^n b_{m,1}\bigg) = 0,
$$
which is also a tautology because of \eqref{eq_zero_sum_coeffs_conj}.  Therefore, the recurrence relation imposes no restriction on $a_1$.  

For $j = 2$, the indices $(j_1,\ldots,j_n)$ in \eqref{eq_Polya_n_recur_conj} are of two types:~either (i) those for which some $j_m = 2$ and all other $j_r = 0$, or (ii) those for which $j_p = j_q = 1$ for some $p < q$ and all other $j_r = 0$.  

In case (i), $\prod_{B \in \nu} a_{j_B} = a_2$ for all $\nu \in \Pi_n$, and the inner sum on the left-hand side of \eqref{eq_Polya_n_recur_conj} reduces to 
$$
a_2 \sum_{\nu \in \Pi_n} (-1)^{|\nu|-1} \prod_{B \in \nu} (|B|-1)!,
$$
which, by \eqref{eq_zero_sum_coeffs_conj}, equals $0$.  

In case (ii), let $\mcalP_n$ be the set of partitions $\nu \in \Pi_n$ such that there exist ``\textit{coupled}'' indices $p$ and $q$, both belonging to a common block in $\nu$, and such that $p < q$ and $j_p = j_q = 1$.  Clearly, $\Pi_n = (\Pi_n \cap \mcalP_n) \cup (\Pi_n \setminus \mcalP_n)$, so 
\begin{align}
\label{eq_Polya_n_reduced_1}
\sum_{j_1+\cdots+j_n=2} \bigg[&\sum_{\nu \in \Pi_n} (-1)^{|\nu|-1} \prod_{B \in \nu} (|B|-1)! \: a_{j_B}\bigg] \prod_{m=1}^n b_{m,j_m} \nonumber \\
&= \sum_{j_1+\cdots+j_n=2} \bigg[\sum_{\nu \in \Pi_n \cap \mcalP_n} (-1)^{|\nu|-1} \prod_{B \in \nu} (|B|-1)! \: a_{j_B} \nonumber \\
& \qquad\qquad\qquad\qquad + \sum_{\nu \in \Pi_n \setminus \mcalP_n} (-1)^{|\nu|-1} \prod_{B \in \nu} (|B|-1)! \: a_{j_B}\bigg] \prod_{m=1}^n b_{m,j_m}.
\end{align}
Since $j_1+\cdots+j_n=2$ then 
$$
\prod_{B \in \nu} a_{j_B} = 
\begin{cases}
a_2, & \textrm{if } \nu \in \Pi_n \cap \mcalP_n \\
a_1^2, & \textrm{if } \nu \in \Pi_n \setminus \mcalP_n
\end{cases}
$$
and also, if $p$ and $q$ are coupled indices then 
$$
\prod_{m=1}^n b_{m,j_m} = b_{p,1} b_{q,1}.
$$
Therefore, for $j=2$, it follows by \eqref{eq_Polya_n_reduced_1} that the recurrence relation \eqref{eq_Polya_n_recur_conj} reduces to 
\begin{multline}
\label{eq_Polya_n_reduced_2}
a_2 \bigg[\sum_{\nu \in \Pi_n \cap \mcalP_n} (-1)^{|\nu|-1} \prod_{B \in \nu} (|B|-1)! \\
+ a_1^2 \sum_{\nu \in \Pi_n \setminus \mcalP_n} (-1)^{|\nu|-1} \prod_{B \in \nu} (|B|-1)!\bigg] \bigg(\sum_{1 \le p < q \le n} b_{p,1} b_{q,1}\bigg) = 0.
\end{multline}
We claim that 
\begin{equation}
\label{eq_Polya_n_reduced_3}
\sum_{\nu \in \Pi_n \cap \mcalP_n} (-1)^{|\nu|-1} \prod_{B \in \nu} (|B|-1)! = (n-2)!.
\end{equation}
To see this we proceed as at \eqref{eq_signless_Stirling}, \textit{infra}, deducing that the sum in \eqref{eq_Polya_n_reduced_3} is a weighted sum over all permutations $\sigma \in \mathcal{S}_n$ such that $\sigma$ has exactly $k$ cycles in its cyclic decomposition, there exist distinct indices $p$ and $q$ that belong to the same cycle, the cycles of $\sigma$ form the blocks of $\nu$, and $\sigma$ is assigned the weight $(-1)^{k-1}$.  Consider the cycle which contains $p$ and $q$; if this cycle has size $l$, then choose its remaining $l-2$ elements, arrange its cycle in $(l-1)!$ ways, and permute the remaining $n-l$ elements.  It follows that the signed sum over all permutations of the remaining set of $n-l$ elements is 
$$
\sum_{k=0}^{n-l} (-1)^k \stirlingI{n-l}{k} = 
\begin{cases}
\phantom{-}1, & l = n \\
-1, & l = n-1 \\
\phantom{-}0, & l \le n-2
\end{cases} \ ,
$$
where we have applied the identity \eqref{eq_identity_signless_Stirling}.  Therefore the terms corresponding to $l = n, n-1$ are the only terms that contribute to the sum in \eqref{eq_Polya_n_reduced_3}, so 
$$
\sum_{\nu \in \Pi_n \cap \mcalP_n} (-1)^{|\nu|-1} \prod_{B \in \nu} (|B|-1)! = (n-1)! - (n-2) \cdot (n-2)! = (n-2)!,
$$
and this establishes \eqref{eq_Polya_n_reduced_3}.  

Next, by \eqref{eq_zero_sum_coeffs_conj}, 
\begin{multline*}
\sum_{\nu \in \Pi_n \setminus \mcalP_n} (-1)^{|\nu|-1} \prod_{B \in \nu} (|B|-1)! + \sum_{\nu \in \Pi_n \cap \mcalP_n} (-1)^{|\nu|-1} \prod_{B \in \nu} (|B|-1)! \\
= \sum_{\nu \in \Pi_n} (-1)^{|\nu|-1} \prod_{B \in \nu} (|B|-1)! = 0,
\end{multline*}
so it follows that 
\begin{align}
\label{eq_Polya_n_reduced_4}
\sum_{\nu \in \Pi_n \setminus \mcalP_n} (-1)^{|\nu|-1} \prod_{B \in \nu} (|B|-1)! &= - \sum_{\nu \in \Pi_n \cap \mcalP_n} (-1)^{|\nu|-1} \prod_{B \in \nu} (|B|-1)! \nonumber \\
&= -(n-2)!.
\end{align}
Applying \eqref{eq_Polya_n_reduced_4} and \eqref{eq_Polya_n_reduced_3}, it follows that the recurrence relation \eqref{eq_Polya_n_reduced_2} reduces to 
$$
(n-2)! (a_2 - a_1^2) \bigg(\sum_{1 \le p < q \le n} b_{p,1} b_{q,1}\bigg) = 0,
$$
which shows that $a_2 = a_1^2$ since $b_{m,j} > 0$ for all $m, j$.  

Now suppose, as inductive hypothesis, that $a_k = a_1^k$ for all $k=0,\ldots,j-1$, for fixed $j \ge 3$.  We write the recurrence relation \eqref{eq_Polya_n_recur_conj} as 
\begin{equation}
\label{eq_Polya_n_recur_conj_2}
\sum_{j_1+\cdots+j_n=j} L_n(j_1,\ldots,j_n) \prod_{m=1}^n b_{m,j_m} = 0,
\end{equation}
where $L_n(j_1,\ldots,j_n)$ is defined in \eqref{eq_L_n_sum}.  Recall, by Lemma \ref{lem_L_n_reduction} that $L_n(j_1,\ldots,j_n) = 0$ if $q = 0,1$, and that $L_n(j_1,\ldots,j_n)$ is symmetric in $j_1,\ldots,j_n$.  Moreover, if $q$ of the indices $j_1,\ldots,j_n$ are non-zero then $q \le \min\{j,n\}$.  Therefore \eqref{eq_Polya_n_recur_conj_2} reduces to 
\begin{align*}
&\sum_{\substack{j_1+j_2=j \\ j_1,j_2 \ge 1}} L_n(j_1,j_2,0,\ldots,0) \bigg(\sum_{1 \le m_1 < m_2 \le n} b_{m_1,j_1} b_{m_2,j_2}\bigg) \\
& \ \ + \sum_{\substack{j_1+j_2+j_3=j \\ j_1,j_2,j_3 \ge 1}} L_n(j_1,j_2,j_3,0,\ldots,0) \bigg(\sum_{1 \le m_1 < m_2 < m_3 \le n} b_{m_1,j_1} b_{m_2,j_2} b_{m_3,j_3}\bigg) \\
& \ \ \ \ + \cdots + \sum_{\substack{j_1+\cdots+j_r=j \\ j_1,\ldots,j_r \ge 1}} L_n(j_1,\ldots,j_r,0,\ldots,0) \bigg(\sum_{1 \le m_1 < \cdots < m_r \le n} b_{m_1,j_1} b_{m_2,j_2} \cdots b_{m_r,j_r}\bigg) = 0,
\end{align*}
where $r = \min\{j,n\}$; more concisely, 
$$
\sum_{q=2}^{\min\{j,n\}} \sum_{\substack{j_1+\cdots+j_q=j \\ j_1,\ldots,j_q \ge 1}} L_n(j_1,\ldots,j_q,0,\ldots,0) \sum_{1 \le m_1 < \cdots < m_q \le n} \prod_{k=1}^q b_{m_k,j_k} = 0.
$$
Now applying the reduction formula \eqref{eq_L_n_reduction_n_q} in Lemma \ref{lem_L_n_reduction}, we have 
\begin{equation}
\label{eq_Polya_n_recur_conj_3}
\sum_{q=2}^{\min\{j,n\}} \frac{(n-2)!}{(q-2)!} \sum_{\substack{j_1+\cdots+j_q=j \\ j_1,\ldots,j_q \ge 1}} L_q(j_1,\ldots,j_q) \sum_{1 \le m_1 < \cdots < m_q \le n} \prod_{k=1}^q b_{m_k,j_k} = 0.
\end{equation}
For $j_1+\cdots+j_q=j$ and $j_1,\ldots,j_q \ge 1$, we have $j_B < j$ for any block $B \in \Pi_q$ such that $|B| \neq q$, i.e., $B \neq \{1,\ldots,q\}$, which is equivalent to $|\nu| \neq 1$.  For each such block $B$, it follows by the inductive hypothesis that $a_{j_B} = a_1^{j_B}$, and then we obtain 
\begin{align*}
L_q(j_1,\ldots,j_q) &= \sum_{\nu \in \Pi_q} (-1)^{|\nu|-1} \prod_{B \in \nu} (|B|-1)! \: a_{j_B} \\
&= (q-1)! \: a_j + \sum_{\substack{\nu \in \Pi_q \\ |\nu| \ge 2}} (-1)^{|\nu|-1} \prod_{B \in \nu} (|B|-1)! \: a_1^{j_B} \\
&= (q-1)! \: a_j + a_1^j \sum_{\substack{\nu \in \Pi_q \\ |\nu| \ge 2}} (-1)^{|\nu|-1} \prod_{B \in \nu} (|B|-1)!,
\end{align*}
since $\prod_{B \in \nu} a_1^{j_B} = a_1^{j_1+\cdots+j_q} = a_1^j$ for $\nu \in \Pi_q$.  By \eqref{eq_zero_sum_coeffs_conj}, 
$$
(q-1)! + \sum_{\substack{\nu \in \Pi_q \\ |\nu| \ge 2}} (-1)^{|\nu|-1} \prod_{B \in \nu} (|B|-1)! \equiv \sum_{\nu \in \Pi_q} (-1)^{|\nu|-1} \prod_{B \in \nu} (|B|-1)! = 0,
$$
hence 
$$
\sum_{\substack{\nu \in \Pi_q \\ |\nu| \ge 2}} (-1)^{|\nu|-1} \prod_{B \in \nu} (|B|-1)! = -(q-1)!.
$$
Therefore under the inductive hypothesis, we have 
$$
L_q(j_1,\ldots,j_q) = (q-1)! \: (a_j - a_1^j),
$$
and by substituting this result in \eqref{eq_Polya_n_recur_conj_3}, we find that the recurrence relation reduces to 
$$
\sum_{q=2}^{\min\{j,n\}} \frac{(n-2)!}{(q-2)!} (q-1)! \sum_{\substack{j_1+\cdots+j_q=j \\ j_1,\ldots,j_q \ge 1}} (a_j - a_1^j) \sum_{1 \le m_1 < \cdots < m_q \le n} \prod_{k=1}^q b_{m_k,j_k} = 0,
$$
which simplifies to 
$$
(a_j - a_1^j) \sum_{q=2}^{\min\{j,n\}} (q-1) \sum_{\substack{j_1+\cdots+j_q=j \\ j_1,\ldots,j_q \ge 1}} \sum_{1 \le m_1 < \cdots < m_q \le n} \prod_{k=1}^q b_{m_k,j_k} = 0.
$$
Since $b_{m,j} > 0$ for all $m, j$ then we obtain $a_j = a_1^j$.  

By induction, we now obtain $a_j = a_1^j$ for all $j=0,1,2,\ldots$, and then it follows that $\psi(t) = \exp(a_1 t^2)$, $t \in \R$, so $X_1$ is normally distributed with mean $0$.  
$\qed$

\section{Characterizations of some contaminated normal distributions}
\label{sec_contaminated_normal}

The next characterization arises from the case in which the coefficient $a_2$ in \eqref{eq_cf_X} satisfies $a_2 = 2a_1^2 - 1$.  In stating this result, we recall from \cite{Mason_Handscomb} that for $j \in \Z$, the $j$th \textit{Chebyshev polynomial of the first kind}, $T_j(\cdot)$, is defined by 
\begin{equation}
\label{eq_Chebyshev_def}
T_j(\cos\theta) = \cos(j\theta), \qquad \theta \in \R.
\end{equation}

Let $\varepsilon$, $Y_1$, and $Y_2$ be mutually independent random variables such that $\varepsilon$ is uniformly distributed on the set $\{0,1\}$, and $Y_1$ and $Y_2$ are centered normal random variables with variances $\sigma_1^2$ and $\sigma_2^2$, respectively.  Then the random variable $X = (1-\varepsilon) Y_1 + \varepsilon Y_2$ is said to have an \textit{$\varepsilon$-contaminated normal distribution with parameters $\sigma_1^2$ and $\sigma_2^2$}.  

In the following result, we will encounter a constant $a_1 \le -1$.  In that case, both $a_1 + (a_1^2-1)^{1/2}$ and $a_1 - (a_1^2-1)^{1/2}$ are \textit{negative} since their product is positive and their sum is negative.  Then we have the following characterization.

\begin{theorem}
\label{thm_Polya_n_Chebyshev}
Let $n$ be an odd integer, $n \ge 3$.  Let $X_1,\ldots,X_n$ be mutually independent, identically distributed random variables whose characteristic function, $\psi(\cdot)$, is analytic.  Let $R_1,\ldots,R_n$ be mutually independent random variables, all of whose moments are finite, and such that $R = (R_1,\ldots,R_n)$ is independent of $X_1,\ldots,X_n$.  
\begin{enumerate}[label=(\roman*), nosep]
\item \label{gen_polya_thm_Chebyshev_i}
Suppose that $X_1$ has an $\varepsilon$-contaminated normal distribution with parameters $\sigma_1^2 = -2\big(a_1 + (a_1^2-1)^{1/2}\big)$ and $\sigma_2^2 = -2\big(a_1 - (a_1^2-1)^{1/2}\big)$, where $a_1 \le -1$.  Then 
\begin{equation}
\label{eq_Polya_R_n_cf_Chebyshev}
\E \sum_{\nu \in \Pi_n} (-1)^{|\nu|-1} (|\nu|-1)! \exp\!\bigg(\i t \sum_{B \in \nu} \|R_B\| X_{\min(B)}\bigg) = 0, 
\qquad t \in \R.
\end{equation}
{\hskip-30pt}Conversely, if \eqref{eq_Polya_R_n_cf_Chebyshev} holds then 
\item \label{gen_polya_thm_Chebyshev_ii}
$\psi(\cdot)$ is symmetric.  
\item \label{gen_polya_thm_Chebyshev_iii}
In \eqref{eq_cf_X}, suppose that $a_j = T_j(a_1)$, $j = 2,\ldots,n-1$.  Then $a_1 \le -1$ and $X_1$ has an $\varepsilon$-contaminated normal distribution with parameters $\sigma_1^2 = -2\big(a_1 + (a_1^2-1)^{1/2}\big)$ and $\sigma_2^2 = -2\big(a_1 - (a_1^2-1)^{1/2}\big)$.  
\end{enumerate}
\end{theorem}

\begin{example}
\label{ex_Chebyshev_n_3}
{\rm
Suppose that $n = 3$.  Let $X_1,X_2,X_3$ be mutually independent, identically distributed random variables, each having a characteristic function $\psi(t)$ that is analytic.  Let $R_1,R_2,R_3$ be mutually independent random variables that are independent of $X_1,X_2,X_3$ and such that all moments of $R_1,R_2,R_3$ are finite.  If 
\begin{align*}
& \E \Big[\exp\!\big(\i t(R_1^2+R_2^2+R_3^2)^{1/2} X_1\big) \\
& \quad\ - \exp\!\big(\i t((R_1^2+R_2^2)^{1/2} X_1 + R_3 X_3)\big) - \exp\!\big(\i t((R_1^2+R_3^2)^{1/2} X_1 + R_2 X_2)\big) \\
& \qquad\quad - \exp\!\big(\i t((R_2^2+R_3^2)^{1/2} X_2 + R_1 X_1)\big) + 2 \exp\!\big(\i t(R_1 X_1 + R_2 X_2 + R_3 X_3)\big)\Big] = 0,
\end{align*}
for all $t \in \R$, then all odd moments of $X_1$ are equal to zero, and $\psi(t) = \phi(t^2)$ where $\phi(\cdot)$ is of the form \eqref{eq_cf_X}.  Further, if $a_2 = 2a_1^2 - 1$ then $a_1 \le -1$ and $X_1$ has an $\varepsilon$-contaminated normal distribution.  
}\end{example}

\subsection{Preliminaries}
\label{subsec_prelim_contaminated_normal}

Setting $x = \cos\theta$ in \eqref{eq_Chebyshev_def}, we obtain 
\begin{equation}
\label{eq_Tj}
T_j(x) = \frac12 \big[\big(x + (x^2-1)^{1/2}\big)^j + \big(x - (x^2-1)^{1/2}\big)^j\big],
\end{equation}
$x \in \R$, $j \in \Z$ \cite{Mason_Handscomb}.  Let $z = x + (x^2-1)^{1/2}$, then $z^{-1} = x - (x^2-1)^{1/2}$ and by  \eqref{eq_Tj}, 
\begin{equation}
\label{eq_Tj_pm}
T_j(x) = \frac12 (z^j + z^{-j}) \equiv T_{-j}(x).
\end{equation}
Suppose that $x > 1$; if $j \ge 0$ then, by expanding the $j$th powers in \eqref{eq_Tj} and cancelling terms, we find that $T_j(x)$ is a real-valued polynomial of degree $j$.  For $j < 0$, it follows from \eqref{eq_Tj_pm} that $T_j(x)$ is a real-valued polynomial of degree $-j$.  Therefore $T_j(x)$, $x \in \R$, is a real-valued polynomial of degree $|j|$.  

Using the representation \eqref{eq_Tj}, it is straightforward to obtain the exponential generating function:~for all $x,t \in \R$, 
\begin{equation}
\label{eq_T_j_gen_function}
\sum_{j=0}^\infty \frac{T_j(x) t^j}{j!} = \frac12 \Big[\exp\!\big((x + (x^2 - 1)^{1/2})t\big) + \exp\!\big((x - (x^2 - 1)^{1/2})t\big)\Big].
\end{equation}

The next result shows that the polynomials $T_j(\cdot)$ satisfy a cumulant-type identity.  Despite the lengthy history of the Chebyshev polynomials this result, as far as we have been able to ascertain, appears to be new.  

\begin{lemma}
\label{lem_zero_sum_Chebyshev}
Let $n$ be an odd integer, $n \ge 3$.  For nonnegative integers $j_1,\ldots,j_n$, define 
\begin{equation}
\label{eq_zero_sum_Chebyshev}
M_n(x) := \sum_{\nu\in\Pi_n} (-1)^{|\nu|-1}(|\nu|-1)! \prod_{B\in\nu} T_{j_B}(x),
\end{equation}
$x \in \R$.  Then $M_n(x) \equiv 0$.  
\end{lemma}

\noindent
\textit{Proof}.  Let the random variable $V$ be uniformly distributed on the set $\{-1,1\}$.  Then $\E \cos(\theta V) = \cos(\theta)$ and $\E \sin(\theta V) = 0$, $\theta \in \R$.  

For $x \in [-1,1]$, let $\theta = \cos^{-1} x$ in \eqref{eq_zero_sum_Chebyshev} and define the random variables $Y_m = \exp(\i j_m \theta V)$, $m=1,\ldots,n$.  Then for any partition $\nu$ and any block $B \in \nu$, it follows by \eqref{eq_Chebyshev_def} that 
$$
\E \!\bigg(\prod_{m \in B} Y_m\bigg) = \E \exp\big(\i j_B \theta V\big) = \cos\big(j_B \theta\big) = T_{j_B}(x).
$$
Therefore 
$$
M_n(x) = \sum_{\nu \in \Pi_n} (-1)^{|\nu|-1} (|\nu|-1)! \prod_{B\in\nu} \E \!\bigg(\prod_{m \in B} Y_m\bigg) = \kappa(Y_1,\ldots,Y_n),
$$
the joint cumulant of $Y_1,\ldots,Y_n$.  Since $V = \pm 1$ then 
$$
Y_m = \cos(j_m\theta V) + \i \sin(j_m\theta V) = \cos(j_m\theta) + \i V \sin(j_m\theta).
$$
Therefore, by the multilinearity property of the joint cumulants, 
\begin{align*}
M_n(x) &= \kappa\big(\!\cos(j_1\theta) + \i V \sin(j_1\theta),\ldots,\cos(j_n\theta) + \i V \sin(j_n\theta)\big) \\
&= \i^n \bigg(\prod_{m=1}^n \sin(j_m\theta)\bigg) \kappa_n(V),
\end{align*}
where $\kappa_n(V)$ is the $n$th cumulant of $V$.  Since $n$ is odd and $V \eqdist -V$ then 
\begin{equation}
\label{eq_kappa_n_V}
\kappa_n(V) = \kappa_n(-V) = (-1)^n \kappa_n(V) = - \kappa_n(V);
\end{equation}
therefore $\kappa_n(V) = 0$.  Hence $M_n(x) = 0$ for $-1 \le x \le 1$, and since $M_n(x)$ is a polynomial then it is identically zero for all $x \in \R$.  
$\qed$

\subsection{Proof of Theorem \ref{thm_Polya_n_Chebyshev}}
\label{subsec_Polya_n_Chebyshev_i}

\noindent
\textbf{Proof of Theorem \ref{thm_Polya_n_Chebyshev}\ref{gen_polya_thm_Chebyshev_i}}: 
Denoting expectation with respect to any random variable $Y$ by $\E_Y$, it follows from by the definition of the $\varepsilon$-contaminated normal distribution that 
\begin{align}
\label{eq_cf_contam_normal}
\psi(t) &= \E_{Y_1,Y_2} \E_\varepsilon \exp\!\big(\i t((1-\varepsilon) Y_1 + \varepsilon Y_2)\big) \nonumber \\
&= \frac12 \E \big[\exp(\i tY_1) + \exp(\i tY_2)\big] \nonumber \\
&= \frac12 \big[\exp\!\big((a_1 + (a_1^2-1)^{1/2}) t^2\big) + \exp\!\big((a_1 - (a_1^2-1)^{1/2}) t^2\big)\big].
\end{align}
Let $t\in\mathbb R$ and let $\tau = (\tau_1,\ldots,\tau_n)$ be a constant vector.  Also, let $V$ be uniformly distributed on $\{-1,1\}$ and define
$$
U_m = \exp\!\big(t^2 (a_1 + (a_1^2-1)^{1/2} V)\tau_m^2\big),
$$
$m = 1,\ldots,n$.  For any block $B \in \nu$, where $\nu \in\Pi_n$, 
$$
\prod_{m \in B} U_m = \exp\!\big(t^2 (a_1 + (a_1^2-1)^{1/2} V) \|\tau_B\|^2\big);
$$
therefore 
\begin{align*}
\E \bigg(\prod_{m\in B} U_m\bigg) &= \frac12 \Big[\exp\!\big((a_1 + (a_1^2-1)^{1/2}) t^2 \|\tau_B\|^2\big) + 
\exp\!\big((a_1 - (a_1^2-1)^{1/2}) t^2 \|\tau_B\|^2\big)\Big] \\
&= \psi(t \|\tau_B\|),
\end{align*}
and it follows by \eqref{eq_cumulant_sum}, the formula for the joint cumulant, that 
\begin{equation}
\label{eq_cumulant_sum_contam_distn}
\sum_{\nu \in \Pi_n} (-1)^{|\nu|-1} (|\nu|-1)! \prod_{B\in\nu} \psi(t \|\tau_B\|) = \kappa(U_1,\ldots,U_n).
\end{equation}
Since $V \in \{-1,1\}$ then 
$$
\exp\!\big(t^2 (a_1^2-1)^{1/2} \tau_m^2 V\big) = \cosh\!\big(t^2 (a_1^2-1)^{1/2} \tau_m^2\big) + V \sinh\!\big(t^2 (a_1^2-1)^{1/2} \tau_m^2\big);
$$
hence 
$$
U_m = \exp(t^2 a_1 \tau_m^2) \cosh\!\big(t^2 (a_1^2-1)^{1/2} \tau_m^2\big) + V \exp(t^2 a_1 \tau_m^2) \sinh\!\big(t^2 (a_1^2-1)^{1/2} \tau_m^2\big),
$$
$m = 1,\ldots,n$.  Therefore by the multilinearity property of the joint cumulants, 
$$
\kappa(U_1,\ldots,U_n) = \bigg(\prod_{m=1}^n \exp(t^2 a_1 \tau_m^2) \sinh\!\big(t^2 (a_1^2-1)^{1/2} \tau_m^2\big)\bigg) \kappa_n(V).
$$
As shown at \eqref{eq_kappa_n_V}, $\kappa_n(V) = 0$ for odd $n$, so the sum on the left-hand side of \eqref{eq_cumulant_sum_contam_distn} equals zero.  

For each $\nu \in \Pi_n$, the random variables $X_{\min(B)}$, indexed by the blocks $B \in \nu$, are mutually independent; therefore 
\begin{equation}
\label{eq_X_B_cfs}
\E \bigg[\exp\!\bigg(\i t \sum_{B \in \nu} \|R_B\| X_{\min(B)}\bigg) \bigg| \, R\bigg] = \prod_{B \in \nu} \psi(t\|R_B\|),
\end{equation}
and this conditional expectation is integrable with respect to $R$ since $|\psi(t)| \le 1$ for all $t$.  Now inserting \eqref{eq_X_B_cfs} into the sum in \eqref{eq_cumulant_sum_contam_distn}, we obtain 
\begin{multline*}
\E \bigg[\sum_{\nu \in \Pi_n} (-1)^{|\nu|-1} (|\nu|-1)! \exp\!\bigg(\i t \sum_{B \in \nu} \|R_B\| X_{\min(B)}\bigg)\bigg| \, R = \tau\bigg] \\
= \sum_{\nu \in \Pi_n} (-1)^{|\nu|-1} (|\nu|-1)! \prod_{B \in \nu} \psi(t \|R_B\|) = 0
\end{multline*}
for all $\tau$, and by taking expectations with respect to $R$ we obtain \eqref{eq_Polya_R_n_cf_Chebyshev}.  
$\qed$

\bigskip


\noindent
\textbf{Proof of Theorem \ref{thm_Polya_n_Chebyshev}\ref{gen_polya_thm_Chebyshev_ii}}: 
Here, the proof 
is the same as the proof of Theorem \ref{thm_Polya_n}\ref{gen_polya_thm_ii}.  
$\qed$

\bigskip


\noindent
\textbf{Proof of Theorem \ref{thm_Polya_n_Chebyshev}\ref{gen_polya_thm_Chebyshev_iii}}: 
We proceed analogously to \eqref{eq_Polya_R_n_cf_phi}--\eqref{eq_Polya_n_recur} by substituting the Maclaurin expansion \eqref{eq_cf_X} into \eqref{eq_Polya_R_n_cf_Chebyshev}, comparing the coefficients of equal powers of $t$, and obtaining the recurrence relation 
\begin{equation}
\label{eq_Polya_n_recur_Chebyshev}
\sum_{j_1+\cdots+j_n=j} \bigg[\sum_{\nu \in \Pi_n} (-1)^{|\nu|-1} (|\nu|-1)! \prod_{B \in \nu} a_{j_B}\bigg] \prod_{m=1}^n b_{m,j_m} = 0,
\end{equation}
$j=0,1,2,\ldots$.  

By \eqref{eq_Tj} there holds $T_0(x) = 1$, $T_1(x) = x$, and $T_2(x) = 2x^2 - 1$.  Since we have set $a_2 = 2a_1^2 -1$ then it follows that $a_j = T_j(a_1)$ for $j=0,1,2$.  

By Lemma \ref{lemma_zero_alt_sum}, if any $j_m = 0$ then the inner sum in square brackets in \eqref{eq_Polya_n_recur_Chebyshev} equals $0$.  In particular, for $j \le n-1$, in which case $j_m = 0$ for at least one $m$, \eqref{eq_Polya_n_recur_Chebyshev} becomes a tautology and therefore induces no restrictions on $a_1,\ldots,a_{n-1}$.  Hence we may assign $a_1,\ldots,a_{n-1}$ freely (subject to the requirement that the resulting $\phi(\cdot)$ is a characteristic function). Thus, we set $a_j = T_j(a_1)$ for $j=0,1,2,\ldots,n-1$.  

For $j \ge n$, since the inner sum in \eqref{eq_Polya_n_recur_Chebyshev} equals zero if any $j_m=0$ then the limits on the outer sum restrict further to $j_1,\ldots,j_n \ge 1$.  Also, we decompose the inner sum into the cases in which $|\nu| = 1$, i.e., $\nu = \{1,\ldots,n\}$, and those for which $|\nu| \ge 2$.  If $|\nu| = 1$ then the corresponding summand is $a_{j_{\{1,\ldots,n\}}} \equiv a_{j_1+\cdots+j_n} = a_j$, so \eqref{eq_Polya_n_recur_Chebyshev} reduces for $j \ge n$ to 
\begin{equation}
\label{eq_Polya_n_recur_reduced}
\sum_{\substack{j_1+\cdots+j_n = j \\ j_1,\ldots,j_n \ge 1}} \bigg[a_j + \sum_{\substack{\nu \in \Pi_n \\ |\nu| \ge 2}} (-1)^{|\nu|-1} (|\nu|-1)! \prod_{B \in \nu} a_{j_B}\bigg] \prod_{m=1}^n b_{m,j_m} = 0.
\end{equation}
We now claim that, for all $j=0,1,2,3,\ldots$, 
\begin{equation}
\label{eq_a_j_Chebyshev}
a_j = T_j(a_1).
\end{equation}
Suppose, as inductive hypothesis, that $a_k = T_k(a_1)$ for all $k=0,1,2,\ldots,j-1$, where $j \ge n$.  Since $j_1+\cdots+j_n = j$ and $j_1,\ldots,j_n \ge 1$ then $j_B \le j-1$ for all blocks $B \in \nu$ with $|\nu| \ge 2$.  By the inductive assumption, $a_{j_B} = T_{j_B}(a_1)$ whenever $B \in \nu$ with $|\nu| \ge 2$, and then we obtain 
$$
\prod_{B \in \nu} a_{j_B} = \prod_{B \in \nu} T_{j_B}(a_1).
$$
Hence 
$$
\sum_{\substack{\nu \in \Pi_n \\ |\nu| \ge 2}} (-1)^{|\nu|-1} (|\nu|-1)! \prod_{B \in \nu} a_{j_B} = \sum_{\substack{\nu \in \Pi_n \\ |\nu| \ge 2}} (-1)^{|\nu|-1} (|\nu|-1)! \prod_{B \in \nu} T_{j_B}(a_1) = - T_j(a_1)
$$
since, by Lemma \ref{lem_zero_sum_Chebyshev}, 
$$
T_j(a_1) + \sum_{\substack{\nu \in \Pi_n \\ |\nu| \ge 2}} (-1)^{|\nu|-1} (|\nu|-1)! \prod_{B \in \nu} T_{j_B}(a_1) \equiv \sum_{\nu \in \Pi_n} (-1)^{|\nu|-1} (|\nu|-1)! \prod_{B \in \nu} T_{j_B}(a_1) = 0.
$$
Therefore, under the inductive hypothesis, \eqref{eq_Polya_n_recur_reduced} becomes 
$$
\big[a_j - T_j(a_1)\big] \sum_{\substack{j_1+\cdots+j_n=j \\ j_1,\ldots,j_n \ge 1}} \prod_{m=1}^n b_{m,j_m} = 0,
$$
and since $b_{m,j} > 0$ for all $m$ and $j$ then it follows that $a_j = T_j(a_1)$.  Therefore, by induction, \eqref{eq_a_j_Chebyshev} holds for all $j=0,1,2,3,\ldots$.  

Next, we identify the characteristic function $\phi(t^2)$ and the corresponding distribution of $X_1$.  On applying the exponential generating function \eqref{eq_T_j_gen_function}, we obtain 
\begin{align}
\label{eq_cf_contaminated}
\phi(t^2) &= \sum_{j=0}^\infty \frac{T_j(a_1) t^{2j}}{j!} \nonumber \\
&= \frac12 \big[\exp\!\big((a_1 + (a_1^2 - 1)^{1/2})t^2\big) + \exp\!\big((a_1 - (a_1^2 - 1)^{1/2})t^2\big)\big].
\end{align}
Recall that $a_1 < 0$, so the remaining values of $a_1$ fall into two cases: (a) $-1 < a_1 < 0$, and (b) $a_1 \le -1$.  

(a) Suppose that $-1 < a_1 < 0$.  Then $(a_1^2 - 1)^{1/2}$ is an imaginary complex number, and therefore $(a_1^2 - 1)^{1/2} = (1-a_1^2)^{1/2} \,\i$.  Substituting this result into \eqref{eq_cf_contaminated}, we obtain 
$$
\phi(t^2) 
= \exp(a_1 t^2) \cos\!\big((1-a_1^2)^{1/2}t^2\big).
$$
Since $a_1 < 0$ then $\phi(t^2) \in L^1(\R)$, so its inverse Fourier transform exists.  Applying a formula from Gradshteyn and Ryzhik \cite[Eq.~3.923(2), p.~494]{Gradshteyn_Ryzhik}, we find that the inverse Fourier transform of $\phi(t^2)$ is 
\begin{equation}
\label{eq_inverse_F_transf}
\frac{1}{2\pi} \int_{-\infty}^\infty \exp(-\i tx) \phi(t^2) \dd t 
= \frac{1}{2\pi^{1/2}} \exp\left(\frac{a_1 x^2}{4}\right) 
\cos\!\bigg(\frac{\theta}{2} - \frac{(1-a_1^2)^{1/2} x^2}{4}\bigg),
\end{equation}
where $\theta = \tan^{-1}\big(-(1-a_1^2)^{1/2}/a_1\big) = \cos^{-1} (-a_1)$.  Since 
$$
\frac{\theta}{2} - \frac{(1-a_1^2)^{1/2} x^2}{4} \to -\infty
$$
as $|x| \to \infty$ then the the right-hand side of \eqref{eq_inverse_F_transf} attains negative values for sufficiently large $|x|$.  Therefore $\phi(t^2)$ in this case is not a characteristic function.  

(b) Suppose that $a_1 \le -1$.  Since both $a_1 - (a_1^2 - 1)^{1/2}$ and $a_1 + (a_1^2 - 1)^{1/2}$ are negative then each function $\exp\!\big((a_1 - (a_1^2 - 1)^{1/2})t^2\big)$ and $\exp\!\big((a_1 + (a_1^2 - 1)^{1/2})t^2\big)$ is a characteristic function, so $\phi(t^2)$ in \eqref{eq_cf_contaminated} is a characteristic function.  Finally, as shown earlier in \eqref{eq_cf_contam_normal}, $X_1$, the corresponding random variable, has a $\varepsilon$-contaminated normal distribution with parameters $-2(a_1 + (a_1^2 - 1)^{1/2})$ and $-2(a_1 - (a_1^2 - 1)^{1/2})$.  
$\qed$

\section{Characterizations of the exponential distribution}
\label{sec_exponential_distn}

Let $X_1,\ldots,X_n$ be mutually independent, identically distributed, nonnegative random variables.  It is well-known that $X_1 \eqdist n \min(X_1,\ldots,X_n)$ if $X_1$ has an exponential distribution; indeed this result has been traced by Marshall and Olkin \cite{Marshall_Olkin} to as far back as Fisher and Tippett \cite{Fisher_Tippett}.  

Let $\theta_1,\ldots,\theta_n$ be positive constants.  The characterization that we shall study states that $X_1$ has an exponential distribution if and only if 
\begin{equation}
\label{eq_contam_exponential_charac}
\frac{X_1}{\theta_1+\cdots+\theta_n} \eqdist \min\left\{\frac{X_1}{\theta_1},\ldots,\frac{X_n}{\theta_n}\right\}.
\end{equation}
For the case in which $n=2$, this characterization is provided by Bryc \cite[Section 3.4]{Bryc}, and the details given there extend easily to cover the case of general $n$.  

To generalize this characterization, we adopt a philosophy similar to that developed in the earlier sections.  In particular, we will replace the constants $\theta_1,\ldots,\theta_n$ by positive, mutually independent random variables $R_1,\ldots,R_n$ that are independent of $X_1,\ldots,X_n$, and we generalize \eqref{eq_contam_exponential_charac} to include terms that correspond to each partition $\nu \in \Pi_n$.  First, we provide the regularity conditions that are required in the sequel.

\subsection{Preliminaries}
\label{subsec_exponential_distn}

Let $\Fbar(t) = \P(X_1 \ge t)$, $t \ge 0$, denote the \textit{survival function} of $X_1$.  We assume that $\Fbar$ has a Maclaurin expansion, 
\begin{equation}
\label{eq_Fbar_expansion}
\Fbar(t) = \sum_{j=0}^\infty \frac{a_j}{j!} t^j,
\end{equation}
valid for all $t \ge 0$ where, necessarily, $a_0 = 1$ and $a_1 \le 0$.  We also assume that there exist positive constants $c_1, c_2$ such that, for all $j=0,1,2,\ldots$, 
\begin{equation}
\label{eq_a_j_bound}
|a_j| \le c_1 c_2^j.
\end{equation}
As for the random variables $R_1,\ldots,R_n$, we assume that there exists $\eta > 0$ such that 
\begin{equation}
\label{eq_R_bound}
\E \exp(\eta R_m) < \infty, \qquad m=1,\ldots,n.
\end{equation}
It follows from \eqref{eq_a_j_bound}, \eqref{eq_R_bound}, and the mutual independence of $R_1,\ldots,R_n$ that, for $\nu \in \Pi_n$ and $|t| < \eta/c_2$, 
\begin{align}
\label{eq_E_tail_bound}
\bigg|\E \prod_{B \in \nu} \bigg(\sum_{j=0}^\infty \frac{a_j}{j!} t^j R_B^j\bigg) \bigg| &\le \E \prod_{B \in \nu} \bigg(\sum_{j=0}^\infty \frac{|a_j|}{j!} |t|^j R_B^j\bigg) \nonumber \\
&\le \E \prod_{B \in \nu} \bigg(c_1 \sum_{j=0}^\infty \frac{\eta^j R_B^j}{j!}\bigg) \nonumber \\
&= c_1^{|\nu|} \E \prod_{B \in \nu} \exp(\eta R_B) 
= c_1^{|\nu|} \prod_{m=1}^n \E \exp(\eta R_m) < \infty. 
\end{align}


\begin{theorem}
\label{thm_exponential_character}
Let $X_1,\ldots,X_n$ be nonnegative, identically distributed random variables, each having the survival function $\Fbar(\cdot)$ given by \eqref{eq_Fbar_expansion} and whose coefficients satisfy \eqref{eq_a_j_bound}.  Let $R = (R_1,\ldots,R_n)$ be a vector of nonnegative random variables, each satisfying \eqref{eq_R_bound}, and suppose that $R_1,\ldots,R_n,X_1,\ldots,X_n$ are mutually independent.  If $X_1$ has an exponential distribution then
\begin{equation}
\label{eq_exponential_charac}
\sum_{\nu \in \Pi_n} (-1)^{|\nu|-1} (|\nu|-1)! \: \P\!\left(\min\!\left\{\frac{X_{\min(B)}}{R_B}: B \in \nu\right\} \ge t\right) = 0, \qquad t \ge 0.
\end{equation}

Conversely, suppose that \eqref{eq_exponential_charac} holds and that, in \eqref{eq_Fbar_expansion}, $a_j = a_1^j$ for $j = 2,\ldots,n-1$.  Then $X_1$ has an exponential distribution.  
\end{theorem}

\begin{example}
{\rm
Suppose that $n=3$.  Then \eqref{eq_exponential_charac} reduces to 
\begin{multline*}
\P\bigg(\frac{X_1}{R_1+R_2+R_3} \ge t\bigg) - \P\left(\min\left\{\frac{X_1}{R_1+R_2},\frac{X_3}{R_3}\right\} \ge t\right) \\
- \P\left(\min\left\{\frac{X_1}{R_1+R_3},\frac{X_2}{R_2}\right\} \ge t\right) - \P\left(\min\left\{\frac{X_1}{R_1},\frac{X_2}{R_2+R_3}\right\} \ge t\right) \\
+ 2 \, \P\left(\min\left\{\frac{X_1}{R_1},\frac{X_2}{R_2},\frac{X_3}{R_3}\right\} \ge t\right) \equiv 0.
\end{multline*}
It is simple to verify that if $X_1$ has an exponential distribution then each probability term, $\P(\cdot)$, above equals $\E\exp(-\lambda(R_1+R_2+R_3)t)$, and therefore the sum equals zero.  

We also note that if in \eqref{eq_exponential_charac} the coefficients $(-1)^{|\nu|-1} (|\nu|-1)!$ is replaced by the constants $c_\nu = (-1)^{|\nu|-1} \prod_{B \in \nu} (|B|-1)!$ that appeared in Section \ref{sec_modified_Polya}, then we will obtain a characterization similar to Theorem \ref{thm_Polya_n_conjugate}.  
}\end{example}

\subsection{Proof of Theorem \ref{thm_exponential_character}}
\label{subsec_proof_thm_exponential_character}

Suppose that $X_1$ has an exponential distribution with parameter $\lambda > 0$.  Since the random variables $X_1,\ldots,X_n,R_1,\ldots,R_n$ are mutually independent and $X_1,\ldots,X_n$ are identically distributed then 
\begin{align*}
\P\!\left(\frac{X_{\min(B)}}{R_B} \ge t\right) &= \P(X_{\min(B)} \ge t R_B) \\
&= \E \P(X_1 \ge t R_B \big| R_B) = \E \exp(- \lambda R_B t).
\end{align*}
Since distinct blocks are disjoint then the random vectors $\{R_B: B \in \nu\}$ are mutually independent.  Therefore 
\begin{align*}
\P\!\left(\min\!\left\{\frac{X_{\min(B)}}{R_B}: B \in \nu\right\} \ge t\right) &= \prod_{B \in \nu} \P\!\left(\frac{X_{\min(B)}}{R_B} \ge t\right) \\
&= \E \prod_{B \in \nu} \exp(- \lambda R_B t) = \E \exp\!\bigg(-\lambda t \sum_{m=1}^n R_m\bigg).
\end{align*}
Substituting this result into the left-hand side of \eqref{eq_exponential_charac} and then applying \eqref{eq_zero_sum_coeffs}, we obtain $0$, the right-hand side of \eqref{eq_exponential_charac}.  

Conversely, suppose that \eqref{eq_exponential_charac} holds.  
Since $X_1,\ldots,X_n,R_1,\ldots,R_n$ are mutually independent and the blocks $B \in \nu$ are disjoint, then the random variables $\{X_{\min(B)}, R_B: B \in \nu\}$, are mutually independent.  Therefore 
\begin{align*}
\P\!\left(\min\!\left\{\frac{X_{\min(B)}}{R_B}: B \in \nu\right\} \ge t\right) 
&= \prod_{B \in \nu} \P\!\left(\frac{X_{\min(B)}}{R_B} \ge t\right) \\
&= \E \prod_{B \in \nu} \P(X_{\min(B)} \ge t R_B | R_B)
= \E \prod_{B \in \nu} \Fbar(R_B t),
\end{align*}
so the left-hand side of \eqref{eq_exponential_charac} equals 
\begin{equation}
\label{eq_exponential_charac_LHS}
\E \sum_{\nu \in \Pi_n} (-1)^{|\nu|-1} (|\nu|-1)! \prod_{B \in \nu} \Fbar(R_B t).
\end{equation}
Next we apply the expansion \eqref{eq_Fbar_expansion} to each term $\Fbar(\cdot)$.  By \eqref{eq_E_tail_bound}, 
$$
\bigg|\E \prod_{B \in \nu} \Fbar(R_B t)\bigg| \equiv \bigg|\E \prod_{B \in \nu} \bigg(\sum_{j=0}^\infty \frac{a_j}{j!} t^j R_B^j\bigg) \bigg| < \infty,
$$
for $0 \le t < \eta/c_2$, and by dominated convergence and Tonelli's theorem we can take term-by-term expectations of the resulting multiple series.  

Therefore \eqref{eq_exponential_charac_LHS}, with $t$ replaced by $z \in C$, is an analytic function on $|z| < \eta/A$.  Further, by assumption, \eqref{eq_exponential_charac_LHS} is identically $0$ for all $t \ge 0$, so the uniqueness theorem for analytic function implies that all its Taylor coefficients at $0$ equal $0$.  

Proceeding as in the proof of Theorem \ref{thm_Polya_n}\ref{gen_polya_thm_iii}, specifically at \eqref{eq_Polya_R_n_cf_phi}--\eqref{eq_Polya_n_recur}, and applying \eqref{eq_E_tail_bound} to reverse the expectation and summations, we obtain from \eqref{eq_exponential_charac} and \eqref{eq_exponential_charac_LHS},
$$
\sum_{j=0}^\infty t^j \sum_{j_1+\cdots+j_n = j} \bigg[\sum_{\nu \in \Pi_n} (-1)^{|\nu|-1} (|\nu|-1)! \prod_{B \in \nu} a_{j_B}\bigg] \prod_{m=1}^n \frac{\E(R_m^{j_m})}{j_m!} = 0,
$$
$t \ge 0$, $j=0,1,2,\ldots$.  Using the notation $b_{m,j} := \E(R_m^{j_m})/j_m!$ and comparing the coefficients of $t^j$, we obtain 
$$
\sum_{j_1+\cdots+j_n = j} \bigg[\sum_{\nu \in \Pi_n} (-1)^{|\nu|-1} (|\nu|-1)! \prod_{B \in \nu} a_{j_B}\bigg] \prod_{m=1}^n b_{m,j} = 0,
$$
$j=0,1,2,\ldots$, which is precisely the recurrence relation \eqref{eq_Polya_n_recur}.  Since $a_j = a_1^j$ for all $j \le n-1$, then by following the earlier approach, we obtain $a_j = a_1^j$ for all $j=0,1,2,\ldots$, and therefore $\Fbar(t) = \exp(a_1 t)$, $t \ge 0$.  

Since $X_1$ is non-trivial then $\Fbar(t) \to 0$ as $t \to \infty$; therefore $a_1 < 0$, and $X_1$ has an exponential distribution with parameter $-a_1$.  
$\qed$

\section{Concluding remarks and open problems}
\label{sec_conclusions}

The results obtained in this article raise numerous issues, and we now describe several of them as possibilities for further research.  

The role of M\"obius inversion in our work raises the issue of whether the results obtained here can be extended to more abstract forms of the classical cumulants given by, e.g., Di Nardo \cite{DiNardo} and Doubilet, Rota, and Stanley \cite{Doubilet_etal}.  Moreover, the connection between abstract M\"obius inversion, Euler characteristics, and simplicial complexes, given in \cite[Section 3.8]{Stanley}, raises the issue of whether cumulant-based characterizations can be related to the work of Naiman and Wynn \cite{Naiman_Wynn}.  

It is an open problem to extend the results obtained here to the stable, non-Gaussian distributions.  Such results would generalize the deep theory developed by, e.g., Shimizu and Davies \cite{Shimizu_Davies}, and would also require weakening the regularity assumption that the characteristic function $\psi(\cdot)$ is analytic.  

In considering extensions to other non-Gaussian distributions circular data, we recall that L\'evy \cite{Levy} proved an analog of P\'olya's theorem for distributions on the unit circle.  This raises the issue of whether L\'evy's theorem can be generalized along the lines of the results obtained here.  

We note that Lukacs \cite[Chapter 8]{Lukacs_75} characterized the Wiener process by means of methods that are reminiscent to us of some approaches to P\'olya's theorem.  Thus we raise the problem of extending Lukacs' results, perhaps along the lines given here.  

Finally, stability theorems for many classical characterizations have now been obtained; see, e.g., \cite{Richards81,YY2}.  In the spirit of that literature, we ask for stability versions of the characterizations obtained here.  

\medskip

\noindent
\textbf{Acknowledgments}.  We are grateful to an anonymous reviewer of our article \cite{Bagyan_Richards}, whose comments motivated us to investigate the topics considered here.  

We are also deeply grateful to Tilmann Gneiting and his colleagues at the Institute for Stochastics and the Institute for Statistics at the Karlsruhe Institute of Technology, and the Heidelberg Institute for Theoretical Studies, for the \textit{herzlich und liebensw{\"u}rdig} 
hospitality given to D.R.~in July, 2026 when this research was underway.  

The authors used ChatGPT and Google Scholar to search for references to the literature on set-partitions and on characterizations of probability distributions, and to proofread the article.

\phantomsection

\addcontentsline{toc}{section}{References}

\bibliographystyle{ims}
\begin{thebibliography}{00}

\bibitem{Avksent'ev}
Avksent'ev, D. I.~(1982).  Equidistribution of random linear forms. (Russian) \textit{Mat. Zametki}, \textbf{31}, 947--948 (English translation: \textit{Math. Notes Acad. Sci. USSR}, \textbf{31}, 477--478).

\bibitem{Bagyan_Richards}
Bagyan, A., and Richards, D.~(2026).  Random linear modulation with spherically symmetric modulators. \textit{Sankhy\=a A}, \textbf{88}, 712--763.  

\bibitem{Bobkov_etal}
Bobkov, S., Chistyakov, G., and G\"otze, F.~(2023). \textsl{Concentration and Gaussian Approximation for Randomized Sums}. Springer, New York

\bibitem{Bryc}
Bryc, W.~(1995). \textsl{Normal Distribution: Characterization with Applications}.  Lecture Notes in Statistics, vol. 100, Springer, Berlin.

\bibitem{Chow_Teicher}
Chow, Y. S., and Teicher, H. (1988).  \textsl{Probability Theory: Independence, Interchangeability, Martingales}, second edition.  Springer, New York.

\bibitem{Comtet}
Comtet, L.~(1974).  \textsl{Advanced Combinatorics: The Art of Finite and Infinite Expansions}, revised edition.  D. Reidel, Dordrecht, Holland.

\bibitem{DiNardo}
Di Nardo, E.~(2015).  Symbolic calculus in mathematical statistics: A review.  \textit{S\'eminaire Lotharingien de Combinatoire}, \textbf{67}, Article B67a.

\bibitem{Doubilet_etal}
Doubilet, P., Rota, G.-C., and Stanley R.~(1972).  On the foundations of combinatorial theory (VI): The idea of generating function.  In: \textsl{Sixth Berkeley Symposium on Mathematical Statistics and Probability, Vol. II: Probability Theory}, pp.~267--318 (eds.: L. Le Cam, J. Neyman, and E. L. Scott).  University of California Press, Berkeley.

\bibitem{Fisher_Tippett}
Fisher, R. A., and Tippett, L. H. C. (1928). Limiting forms of the frequency distributions of largest or smallest member of a sample.  \textit{Proc. Cambridge Philos. Soc.}, \textbf{24}, 180--190.

\bibitem{Ghurye_Olkin}
Ghurye, S. G., and Olkin, I.~(1973). Identically distributed linear forms and the normal distribution. \textit{Adv. Appl. Probab.}, \textbf{5}, 138--152.

\bibitem{Gradshteyn_Ryzhik}
Gradshteyn, I. S., and Ryzhik, I. W.~(2007). \textsl{Tables of Integrals, Series, and Products}, seventh edition (A. Jeffrey and D. Zwillinger, editors). Academic Press, San Diego.

\bibitem{Gupta_Richards}
Gupta,~R.~D., and Richards, D.~St.~P.~(1987).  Multivariate Liouville distributions.  \textit{J. Multivariate Anal.}, \textbf{23}, 233--256.

\bibitem{Hamedani}
Hamedani, G. G.~(1992). Bivariate and multivariate normal characterizations: a brief survey. \textit{Commun. Statist.--Theory \& Methods}, \textbf{21}, 2665--2688.

\bibitem{Hornich}
Hornich, H.~(1942). Eine ungleichung f\"ur vektorl\"angen.  \textit{Math. Zeitschrift}, \textbf{48}, 268--274.

\bibitem{Kagan_etal}
Kagan, A. M., Linnik, Yu. V., and Rao, C. R.~(1973).  \textsl{Characterization Problems of Mathematical Statistics}. Wiley, New York.

\bibitem{Kakosyan_etal}
Kakosyan, A. V., Klebanov, L. B., Melamed, J. A.~(1984). \textsl{Characterization of Distributions by the Method of Intensively Monotone Operators}.  Springer, New York.


\bibitem{Klebanov}
Klebanov, L. B.~(2023). Characterizations of probability distributions by the properties of linear forms with random coefficients. \textit{J. Math. Sciences}, \textbf{273}, 778--786.

\bibitem{Levy}
L\'evy, P.~(1939).  L'addition des variables al\'eatoires d\'efinies sur une circonf\'erence.  \textit{Bull. Soc. Math. France}, \textbf{67}, 1--41.


\bibitem{Lukacs_75}
Lukacs, E. (1975).  \textsl{Stochastic Convergence}, second edition.  Academic Press, New York.

\bibitem{Marshall_Olkin}
Marshall, A. W., and Olkin, I.~(1991).  Functional equations for multivariate exponential distributions.  \textit{J. Multivariate Anal.}, \textbf{39}, 209--215.

\bibitem{Mason_Handscomb}
Mason, J. C., and Handscomb, D. C.~(2002). \textsl{Chebyshev Polynomials}. Chapman \& Hall/CRC, Boca Raton, FL.

\bibitem{McCullagh}
McCullagh, P.~(2018).  \textsl{Tensor Methods in Statistics}, second edition.  Dover Publications, Mineola, NY.

\bibitem{Naiman_Wynn}
Naiman, D. Q., Wynn, H. P.~(1992).  Inclusion-exclusion-Bonferroni identities and inequalities for discrete tube-like problems via Euler characteristics. \textit{Ann. Statist.}, \textbf{20}, 43--76.

\bibitem{Polya}
P\'olya, G.~(1923).  Herleitung des Gau{\ss}schen fehlergesetzes aus einer funktionsalsgleichung.  \textit{Math. Zeitschrift}, \textbf{18}, 96--108.

\bibitem{Ressel}
Ressel, P.  (2015).  The Hornich-Hlawka inequality and Bernstein functions.  \textit{J. Math. Inequal.}, \textbf{9}, 883--888.

\bibitem{Richards81}
Richards, D. St. P.~(1981).  Stability theorems for some characterizations of the exponential distribution.  \textit{Ann. Inst. Statist. Math.}, \textbf{33}, 199--204.

\bibitem{Richards}
Richards, D. St. P.~(2004).  Algebraic methods toward higher-order probability inequalities, II.  \textit{Ann. Probab.}, \textbf{32}, 1509--1544.

\bibitem{Sahi}
Sahi, S.~(2008).  Higher correlation inequalities.  \textit{Combinatorica}, \textbf{28}, 209--227.

\bibitem{Shimizu_Davies}
Shimizu, R., and Davies, L. (1981).  General characterization theorems for the Weibull and the stable distributions.  \textit{Sankhy\=a A}, \textbf{43}, 282--310.

\bibitem{Stanley}
Stanley, R. P.~(2012).  \textsl{Enumerative Combinatorics}, Volume 1, second edition.  Cambridge University Press, New York.
 
\bibitem{YY2}
Yanushkevichius, R., and Yanushkevichiene, O.~(2007). Stability of a characterization by the identical distribution of linear forms.  \textit{Statistics}, \textbf{41}, 345--362.

\bibitem{Zaitsev_etal}
Zaitsev, A.~Y., Kagan, A.~M., and Nikitin, Y.~Y.~(2019). Toward the history of the St. Petersburg school of probability and statistics.~IV.~Characterization of distributions and limit theorems in statistics.  \textit{Vestnik St.Petersb. Univ. Math.}, \textbf{52}, 36--53.

\end{thebibliography}

\end{document}